\documentclass[journal]{IEEEtran}

\newcommand{\beq}{\begin{equation}}
\newcommand{\eeq}{\end{equation}}
\newcommand{\bsq}{\begin{subequations}}
	\newcommand{\esq}{\end{subequations}}
\newcommand{\bq}{\begin{eqnarray}}
\newcommand{\eq}{\end{eqnarray}}
\newcommand{\bqn}{\begin{eqnarray*}}
	\newcommand{\eqn}{\end{eqnarray*}}
\usepackage{bbm}
\usepackage[pdftex]{graphicx}
\usepackage{enumerate}
\usepackage{enumitem} 

\usepackage{mathptmx}       
\DeclareMathAlphabet{\mathcal}{OMS}{cmsy}{m}{n}
\usepackage{helvet}         
\usepackage{courier}        
\usepackage{type1cm}        
\usepackage{makeidx}         
\usepackage{graphicx}        

\usepackage{multicol}        
\usepackage[bottom]{footmisc}
\usepackage{ifthen}
\usepackage{subfigure}
\usepackage[usenames,dvipsnames]{color}

\makeindex             

\usepackage{graphics} 
\usepackage{graphicx} 
\usepackage{epsfig} 
\usepackage{epstopdf}
\usepackage{mathptmx} 
\usepackage{times} 
\usepackage[cmex10]{amsmath}

\usepackage{mathtools}  
\usepackage{amssymb}  
\usepackage{amsthm}

\usepackage{amsfonts}
\usepackage{cite}
\usepackage{verbatim}
\usepackage{comment}
\usepackage{algorithmic}
\usepackage{multirow}
\usepackage{cite}
\usepackage{stfloats}
\usepackage{supertabular}
\usepackage{longtable}

\DeclareGraphicsExtensions{.pdf,.png,.jpg,.eps,.ps}
\renewcommand{\arraystretch}{1.2}
\usepackage{arydshln}
\usepackage{multicol}
\usepackage{colortbl}
\theoremstyle{definition}

\newtheorem{lemma}{Lemma}

\newtheorem{proposition}{Proposition}

\theoremstyle{definition}
\newtheorem{definition}{Definition}

\ifCLASSINFOpdf
\else
\fi

\usepackage{comment}
\usepackage{ifthen}
\newboolean{showcomments}
\setboolean{showcomments}{true}

\usepackage{amsmath}
\allowdisplaybreaks[4]

\usepackage[ruled]{algorithm2e}

\begin{document}

%
\title{Approaching Prosumer Social Optimum via Energy Sharing with Proof of Convergence}
%
%
%

\author{Yue Chen, Changhong Zhao, Steven H. Low, and Shengwei Mei
}

%
%

\markboth{Journal of \LaTeX\ Class Files,~Vol.~XX, No.~X, Feb.~2019}%
{Shell \MakeLowercase{\textit{et al.}}: Bare Demo of IEEEtran.cls for IEEE Journals}
%



\maketitle

\begin{abstract}
With the advent of prosumers, the traditional centralized operation may become impracticable due to computational burden, privacy concerns, and conflicting interests. In this paper, an energy sharing mechanism is proposed to accommodate prosumers' strategic decision-making on their self-production and demand in the presence of capacity constraints. Under this setting, prosumers play a generalized Nash game. We prove main properties of the game: an equilibrium exists and is partially unique; no prosumer is worse off by energy sharing and the price-of-anarchy is $1-O(1/I)$ where $I$ is the number of prosumers. In particular, the PoA tends to 1 with a growing number of prosumers, meaning that the resulting total cost under the proposed energy sharing approaches social optimum. We prove that the corresponding prosumers' strategies converge to the social optimal solution as well. Finally we propose a bidding process and prove that it converges to the energy sharing equilibrium under mild conditions. Illustrative examples are provided to validate the results. 
\end{abstract}

\begin{IEEEkeywords}
	Energy sharing, generalized Nash equilibrium, prosumer, bidding algorithm, distributed mechanism
\end{IEEEkeywords}

%
\IEEEpeerreviewmaketitle

\section*{Nomenclature}
	\addcontentsline{toc}{section}{Nomenclature}
	\subsection{Indices, Sets, and Functions}
	\begin{IEEEdescription}[\IEEEusemathlabelsep\IEEEsetlabelwidth{${\underline P _{mn}}$,${\overline P _{mn}}$}]
		\item[$i, \mathcal{I}$] Index and set of prosumers.
		\item[$S_i$] Action sets of prosumer $i$, and $S=\prod_{i \in \mathcal{I}}S_i$.
		\item[$f_i(p_i)$] Cost function of prosumer $i$.
		\item[$u_i(d_i)$] Utility function of prosumer $i$
		\item[$J_i(p_i,d_i)$] Net cost of prosumer $i$, which equals to $f_i(p_i)-u_i(d_i)$; and $J(p,d)=\sum_{i \in \mathcal{I}} J_i(p_i,d_i)$.
		\item[$\Gamma_i(p,d,b)$] Total net cost of prosumer $i$ with sharing, which equals to $f_i(p_i)-u_i(d_i)+\lambda(-a\lambda+b_i)$. 
		\item[PoA$(\mathcal{G})$] Price of anarchy of a game $\mathcal{G}$.
		\item[$\mathcal{Y}_i$] Any $(p_i,d_i) \in \mathcal{Y}_i$ satisfies the capacity constraint.
	\end{IEEEdescription}
	\subsection{Parameters}
	\begin{IEEEdescription}[\IEEEusemathlabelsep\IEEEsetlabelwidth{${\underline P _{mn}}$,${\overline P _{mn}}$}]
		\item[$I$] Number of prosumers.
		\item[$\underline{p}_i,\overline{p}_i$] Lower/upper bound of prosumer $i$'s production.
		\item[$\underline{d}_i,\overline{d}_i$] Lower/upper bound of prosumer $i$'s demand.
		\item[$a$] Energy sharing market sensitivity.
	\end{IEEEdescription}
		\subsection{Decision Variables}
		\begin{IEEEdescription}[\IEEEusemathlabelsep\IEEEsetlabelwidth{${\underline P _{mn}}$,${\overline P _{mn}}$}]
		\item[$p_i$] Production of prosumer $i$.
		\item[$d_i$] Demand of prosumer $i$.
		\item[$\lambda_m$] Dual variable of the power balancing condition.
		\item[$\lambda$] Energy sharing price.
		\item[$q_i$] Amount of energy prosumer $i$ gets from sharing.
		\item[$b_i$] Bid of prosumer $i$ in the energy sharing market.
 		\item[$\tilde{p}_i,\tilde{d}_i$] Optimal strategies under centralized paradigm.
 		\item[$\hat p_i,\hat d_i$] Strategy of prosumer $i$ at sharing equilibrium.
 		\item[$\check{p}_i,\check{d}_i$] Strategy of prosumer $i$ under self-sufficiency.
		\end{IEEEdescription}

\section{Introduction}
%
%
%
%
\IEEEPARstart{I}{n} the US, over 81,000 distributed wind turbines with a cumulative capacity of 1,076 MW had been deployed during 2003-2017 \cite{orrell20162015}. The residential solar photovoltaic (PV) panels had risen from 3,700 MW to 150,000 MW from 2004 to 2014 \cite{agnew2015effect}. Advances in these technologies, together with decline in cost, have encouraged traditional consumers to produce and store energy at home, via distributed energy resources (DERs), electric vehicles, and batteries \cite{parag2016electricity}, turning them into so-called ``prosumers''. They can play a proactive role in energy management. However, a large number of participants, asymmetric information, and conflicting interests also impose great challenges \cite{parag2016electricity}.

\begin{figure}[t]
	\centering
	\includegraphics[width=0.95\columnwidth]{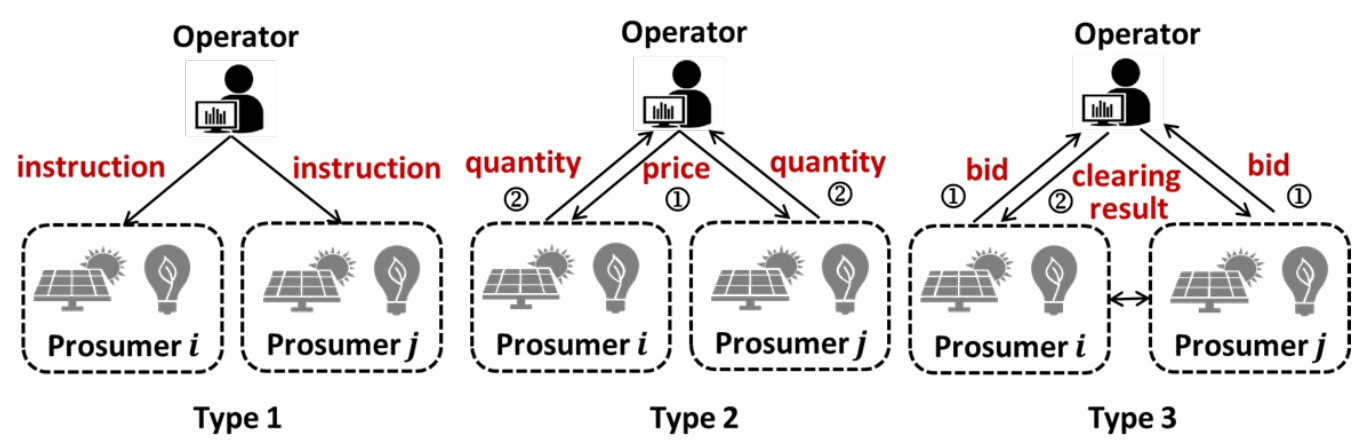}
	\vspace{-2mm}
	\caption{Prosumer management paradigms.}
	\label{fig:paradigm}
\end{figure}

Typically there are three types of prosumer management approaches as shown in Fig. \ref{fig:paradigm} \cite{liu2017energy}. The first one (on the left of Fig. \ref{fig:paradigm}) adopts a centralized operation \cite{zhang2013robust}. The operator of a microgrid or a virtual power plant (VPP) gathers all information and makes a centralized decision, aiming at minimizing the total net costs of all prosumers under management. Then dispatch orders are sent to each prosumer to execute. Since the number of prosumers is increasing rapidly, the traditional centralized approach becomes impracticable both in computational burden and privacy requirements. The second one (in the middle of Fig. \ref{fig:paradigm}) uses a market structure similar to the retail market. The operator announces a price, based on which every prosumer, as a price taker, decides how much energy to consume/produce/buy/sell \cite{liu2017energy2}. It is hard to decide on an effective energy price especially with a large number of prosumers, since private information may be needed and each prosumer's capacity is too small to be observed. Inspired by the concept of ``sharing'' in other sectors, the third approach (on the right of Fig. \ref{fig:paradigm}) has captured increasing attention in recent decades. Here, a prosumer is allowed to exchange/share energy with other prosumers, and they are turning from price-takers to price-makers. This can be done in a peer-to-peer (P2P) structure \cite{hayes2020co} or with the assistance of a platform \cite{chen2019energy}.
As revealed in \cite{chen2018analyzing}, energy sharing can be a promising direction in managing prosumers, since it can achieve a nearly social optimal solution in a distributed manner. Various research projects have been carried out on related issues, such as Piclo \cite{Piclo} in the UK, TransActive Grid \cite{mengelkamp2018designing} in the US, and Enexa \cite{Enexa} in South Australia. The successful operation of energy sharing relies on a well-designed mechanism, and existing research about energy sharing mechanism design can be classified into two categories:

\emph{Cooperative game based approach}. In this approach, first a profit distribution scheme is designed, and then each prosumer chooses its strategy taking into account possible reallocation. The key point here is to design an effective distribution scheme so that all prosumers are willing to collaborate to achieve a certain goal (usually social optimal). Profit distribution schemes for storage sharing were developed under two scenarios \cite{chakraborty2018sharing,han2018incentivizing}. 
Reference \cite{zhou2020distributed} systematically established a pragmatic distributed control and communication mechanism for achieving efficient and resilient coordination of networked microgrids. A mathematical program with equilibrium constraints (MPEC) was used for DER sharing among prosumers, and the coordination surplus was split among the aggregator and prosumers \cite{qi2017sharing}.  
A random sampling method was proposed to estimate the Shapley value of a P2P energy sharing game \cite{han2019estimation}. The trading mechanism with Shapley value was compared with three traditional mechanisms, i.e., bill sharing, mid-market rate, and supply-demand ratio \cite{long2019game}. The cooperative game based sharing can achieve a desired equilibrium with proper distribution rules. However, these rules may require prosumers' private information and thus be difficult to implement.

\emph{Noncooperative game based approach}. This category characterizes the conflicting interests of prosumers, and can be further divided into \emph{bilateral contract based approach} and \emph{auction based approach}. Under the bilateral contract based paradigm, trading offers are posted and handshakes are made. First, each prosumer is registered
as a seller or a buyer. Then, during the trading periods, both the sellers and buyers put forward several offers and try to find the best match. Once a contract is approved by the operator, the corresponding offers are removed \cite{liu2019peer}. This mechanism allows sharing to be performed in an asynchronous manner. A key feature of the bilateral contract based approach is scalability in terms of both the outcomes and the process to reach them. A bilateral contract network with forward and real-time markets was developed in \cite{morstyn2018bilateral}. Reference  \cite{ryu2020real} presented a matching algorithm for microgrid prosumers with minimum risk of mismatch. There are relatively few analytical works based on bilateral contracts, because the matching procedure of bilateral contracts is hard to characterize \cite{ostrovsky2008stability}. 
Under the auction based approach, some or all of the prosumers first bid on energy, and then the market is cleared with the energy sharing prices determined. An evolutionary game was used to model the dynamics of buyers selecting sellers \cite{paudel2018peer}. A Nash bargaining model was adopted in \cite{dutta2014game}  to address the shared charging of electric vehicles (EVs). The energy trading interactions among producers and consumers are modeled as a Stackelberg game, where the producers are leaders and consumers are followers \cite{anoh2019energy}. A shared facilities controller can trade with several residential units at a price in-between the selling and buying prices of the grid, resulting in a win-win game \cite{tushar2014three}. The “Elecbay” platform facilitates the peer-to-peer energy trading was introduced in \cite{zhang2016bidding} with simulation of users’ bidding processes. A double auction mechanism between sellers and buyers was developed in \cite{wang2014game}. In the above works, the role of a prosumer as a buyer or a seller is predetermined and cannot change during the bidding process, so that the sellers can set selling prices, and then the buyers can decide on their demand in response to these prices. However, basically each prosumer could be either a buyer or a seller, depending on other prosumers’ situations that are not known in advance. For example, assume there are three types of prosumers, i.e. High-Cost (HC), Medium-Cost (MC), and Low-Cost (LC). If MC trades with HC, then MC would be a seller while HC would be a buyer. In contrast, if MC trades with LC, then MC would be a buyer while LC would be a seller. Predetermination of market roles will greatly limit the flexibility of prosumers.
To overcome this limit, distributed peer-to-peer energy exchange was modeled as a generalized Nash game in \cite{le2020peer}, which proved that the set of variational equilibria coincides with the set of social optima. In our previous work \cite{chen2019energy}, a generalized demand-bidding approach was proposed for node-level energy sharing, and properties of the Nash equilibrium were proved. This paper extends \cite{chen2019energy} with contributions summarized below.

\subsection*{Contributions}

1) \textbf{Impact}. Over many decades, the traditional power system operation structure has been proven to be quite successful and reliable. Specially, at the demand side, customers are managed by aggregators and usually not price-responsive. With the prevalence of distributed energy resources (DERs), traditional consumers are turning into so-called prosumers, who can trade-off between supply and demand and participate in energy management proactively. In addition, the intermittence and uncertain nature of DERs call for 
stronger capability to deal with real-time energy fluctuation. Exploiting demand-side flexibility to support real-time energy balancing will reduce required generation reserve and save costs. However, the traditional centralized scheme fails to allow a customer to act upon its profit maximizing philosophy, which reduces incentive and restricts demand-side flexibility. Therefore, a new prosumer-oriented approach is desired.

In our previous work \cite{chen2019energy}, an energy sharing mechanism was proposed and several desired properties of the market equilibrium have been proved. We show in \cite{chen2019energy} that the outcome of the proposed energy sharing market approaches that of the centralized operation with an increasing number of prosumers. However, \cite{chen2019energy} merely focuses on the steady-state property (equilibrium) of the energy sharing market, while this paper studies its dynamic property. To be specific, a bidding process in line with each participant’s economic rationality is given in this paper, where each prosumer takes into account the impact of its bid on the sharing price. We reveal that as the number of prosumer grows, the bidding process turns out to have the same form as the Lagrange multiplier-based method for distributed optimization of a centralized problem. This means with more prosumers, not only does its steady-state property converge to the centralized operation, so does its dynamic property. Therefore, many current findings/technologies/theories under centralized operation are likely applicable to the proposed energy sharing market.

2) \textbf{Technical content.} (i) \emph{Model.} Extending \cite{chen2019energy}, in which only power balance constraint was considered with fixed energy demand, this paper incorporates capacity limits and variable demand, which is more flexible and practical. This complicates the analyses in two ways: Firstly, the energy sharing model in this paper can no longer be simplified to a standard Nash game, but indeed is a generalized Nash game whose equilibrium is hard to characterize \cite{facchinei2010generalized}. Secondly, when analyzing each prosumer’s strategic behavior, the complementary slackness conditions associated with the inequality constraints introduce new difficulties. (ii) \emph{Equilibrium.} Main properties of the proposed energy sharing game are proved with three modifications/improvements compared with \cite{chen2019energy}: The generalized Nash equilibrium is partially unique (explained latter), instead of being unique in \cite{chen2019energy}. The energy sharing game achieves a $1-O(1/I)$ price-of-anarchy (PoA, which is less than 1 because net cost is negative in this paper). Besides, as the number of prosumers increases, not only does the total net cost in \cite{chen2019energy}, but also individual prosumers’ strategies, converge to the outcome under the centralized operation, which is a new result not provided in \cite{chen2019energy}. (iii) \emph{Algorithm}. A bidding process is developed for achieving energy sharing in a distributed manner. This paper provides guidance for selecting market sensitivity parameter $a$ so that the bidding process is guaranteed to converge.

3) Our work also differs from \cite{le2020peer}: The prosumers in \cite{le2020peer} have a quadratic cost function and a quadratic utility function, while our design applies to a more general category of strictly convex cost functions and strictly concave utility functions. In \cite{le2020peer}, the prosumers are price-takers in that they decide on their generation or consumption without taking into account the impact of their decisions on the prices. In our paper the prosumers are price-makers, and because of this, the production and demand at generalized Nash equilibrium is the optimal solution of \eqref{eq:central} rather than the social optimal solution in \cite{le2020peer}. Moreover, our work reveals advantageous features of the energy sharing design in Propositions \ref{prop2}-\ref{prop4} and provides a practical bidding process (Algorithm 1), which were not available in \cite{le2020peer}.

\subsection*{Comparison with Relevant Concepts}
\textbf{Energy sharing market \& pool electricity market}. Given that the output of an individual prosumer is too small for it to join the pool electricity market directly, a new approach that allows prosumers to make a profit by exchanging energy with each other is desired. A microgrid connecting those prosumers would be an ideal venue to carry out such an energy exchange, which motivates the energy sharing mechanism proposed in this paper. Specially, in the pool electricity market, a participant is registered in advance as a seller or a buyer, usually a generator as a seller and a load as a buyer. Then the market is cleared by setting a price so that the total supply equals the total demand. However, in an energy sharing market, each prosumer could be either a buyer or a seller, depending on other prosumers’ decisions, so that its role is endogenously determined by the sharing market.

\textbf{The bidding process \& Lagrange multipliers methods}.
The Lagrange multiplier method \cite{shao2016partial} uses the shared multipliers to coordinate different prosumers. The prosumers are ``price takers'' in that they make their generation or consumption decisions without taking into account the impact of their decisions on the multipliers in the next iteration.  This process converges to the socially optimal solution. This model, however, is not applicable to the case where the prosumers are not price takers, but will make strategic decisions that take into account of the impact of their decisions on energy prices. In this case, our proposed bidding process will converge to a generalized Nash equilibrium (GNE) of the energy sharing market (Proposition \ref{prop1}).  Moreover the GNE converges to the social optimal solution as the number of prosumers increases (Proposition \ref{prop3}). We conjecture that our bidding process converges to the Lagrange multipliers method as the number of prosumers increases.


\textbf{Notation}. We use $x:=(x_i, i \in \mathcal{I})^T$ to denote a collection of $x_i$ in a set $\mathcal{I}$. The subscript $-i$ means all components in $\mathcal{I}$ except $i$. The Cartesian product of sets $S_i$ is denoted as $\prod_{i \in \mathcal{I}} S_i$. We use $\dot{f}(.)$ to denote the first derivative of function $f(.)$, and $\ddot{f}(.)$ to denote the second derivative.

\section{Mathematical Formulation}

\subsection{Problem Description}
Consider $I$ prosumers indexed by $i \in \mathcal{I}=\{1,2,...,I\}$ in a standalone microgrid. Assume each prosumer has a distributed generator and a responsive load, whose cost (utility) functions are modeled separately. Specifically, prosumer $i$ produces power $p_i$ at cost $f_i(p_i)$; concurrently, its load consumes power $d_i$ to obtain utility $u_i(d_i)$; function $f_i$ is strictly convex, $u_i$ is strictly concave, and both functions are twice differentiable. Moreover, it is reasonable to assume that $\ddot{f}_i$ and $-\ddot{u}_i$ are uniformly bounded over all $i \in \mathcal{I}$ where both the upper and lower bounds are \emph{strictly} positive and independent from $I$. 

Traditionally, the operator manages all the prosumers in a centralized manner by solving the following problem:
\bsq
\label{eq:opt3}
\begin{align}
\mathop{\min}_{p_i,d_i,\forall i \in \mathcal{I}}~~ & \sum \limits_{i=1}^I [f_i(p_i)-u_i(d_i)] \label{eq:opt3.1}\\
\mbox{s.t.}~~ & \sum \limits_{i=1}^I p_i - \sum \limits_{i=1}^I d_i = 0 \label{eq:opt3.2} :\lambda_m\\
~~ & \underline{p}_i \le p_i \le \overline{p}_i,\forall i \in \mathcal{I} \label{eq:opt3.3}\\
~~ & \underline{d}_i \le d_i \le \overline{d}_i,\forall i \in \mathcal{I} \label{eq:opt3.4}
\end{align}
\esq
where objective \eqref{eq:opt3.1} minimizes the total net cost (cost minus utility) of all the prosumers. Define $J_i(p_i,d_i):=f_i(p_i)-u_i(d_i)$ and $J(p,d)=\sum_{i \in \mathcal{I}} J_i(p_i,d_i)$. 
Constraint \eqref{eq:opt3.2} enforces microgrid-wide power balance (with dual variable $\lambda_m$). Constraints \eqref{eq:opt3.3}-\eqref{eq:opt3.4} impose constant capacity bounds $\underline{p}_i$, $\overline{p}_i$, $\underline{d}_i$, $\overline{d}_i$ on prosumer $i$'s generation and demand. We use problem \eqref{eq:opt3} as the benchmark for subsequent analysis. 

Throughout the paper we assume:

\noindent \textbf{A1}: Problem \eqref{eq:opt3} is feasible.  

Due to strict convexity of objective function, under A1, problem \eqref{eq:opt3} has a unique optimal solution, denoted as $(\tilde{p},\tilde{d})$. 

\textbf{Remark:} 
The dual optimal solution $\tilde{\lambda}_m$ of \eqref{eq:opt3} is unique if there exists at least one prosumer that \emph{strictly} satisfies \eqref{eq:opt3.3} or \eqref{eq:opt3.4} at $(\tilde{p},\tilde{d})$.
A similar discussion about uniqueness of locational marginal price (LMP) can be found in \cite{zhang2014congestion}. The dual optimal $\tilde{\lambda}_m$, known as the ``shadow price'', indicates the increment of total net cost should there be one more unit production-demand mismatch.
We will show that the energy sharing price under the proposed mechanism converges to $\tilde{\lambda}_m$ with a growing number of prosumers.

\subsection{Rationality and Extensions for Assumptions} 
To facilitate theoretical study of fundamental structures of the proposed mechanism, we have made a set of simplifying assumptions, for which the underlying rationality and possible extensions are discussed below. 

(1) \textbf{Cost and utility functions}. We adopt strictly concave utility functions \cite{samadi2012advanced} and strictly convex cost functions \cite{wei2014robust} widely used in power systems. One example for cost function is $f_i(p_i)=\alpha_i^1 p_i^2+\alpha_i^2 p_i$ (with constant parameters $\alpha_i^1,\alpha_i^2>0$), and one for utility function 
is $u_i(d_i)=\beta_i^1d_i^2+\beta_i^2d_i$ (with constant parameters $\beta_i^1<0$, $\beta_i^2>0$).

(2) \textbf{Feasible set}. We allow prosumers to adjust power consumption within capacity limits, following a common simplified demand response model \cite{li2011optimal,samadi2012advanced}. 
The proposed mechanism is compatible with a fixed or precisely predictable demand $d_i^0$ by allowing $\underline d_i=\overline{d}_i=d_i^0$. Moreover, the box constraints \eqref{eq:opt3.3}--\eqref{eq:opt3.4} can be generalized to convex compact sets uniformly bounded over $i \in \mathcal{I}$, with which Propositions 1, 2, 4, and \eqref{eq:PoA} in Proposition 3 still hold. 
In practice, there might be binary variables making the problem nonconvex and thus more challenging. To partially address this concern, our proposed mechanism can be applied 
in concert with appropriate convex relaxation. For instance, the binary variables indicating battery charging or discharging can be converted to complementarity constraints and then tackled  with the exact convex relaxation method in \cite{li2015storage}. 

(3) \textbf{Neglecting network constraints}. Our study is restricted to a residential area or a small microgrid, whose aggregate load only accounts for a small fraction of the total demand at a specific node of a city-sized distribution network. In this case 
it is reasonable to neglect network constraints, as what has been done in \cite{wang2016cooperative,choi2016advanced,mondal2015distributed}.

(4) \textbf{Single time step formulation}. With the prevalence of price-sensitive loads, 
power consumption and price can change dynamically in response to and impacting each other.
Besides, the expanding deployment of DERs, which produce uncertain outputs from geographically dispersed sites, will exert a huge challenge to real-time power balancing. 
In this context, analyzing prosumer behavior and market reliability \emph{in real time} is a crucial topic \cite{chen2020decentralized}. Concerning this topic, we focus on hour-ahead bidding in real-time market, which is commonly modeled as a single time step problem \cite{pei2016optimal}.
Extension to multiple time steps would improve practicality of this work, e.g., to incorporate energy storage; however, it would also tremendously sophisticate notation and presentation while only adding limited value in revealing the fundamental structures and properties of the proposed mechanism. Therefore, we leave this extension for future work.

\subsection{Practical Issues and Requirements}
Although centralized management of prosumers can achieve the lowest total net cost, it encounters two main difficulties in practice: 1) it would be time-consuming when there are a large number of prosumers; 2) information such as $f_i(.)$, $u_i(.)$ is hard to obtain due to privacy concerns of prosumers. To tackle these challenges, a distributed and scalable paradigm is desired, which needs to be:

For process:
1) \emph{Private}. Prosumer privacy is preserved.
2) \emph{Distributed}. Each prosumer makes its own decision based on individual rationality.
3) \emph{Convergent}. The bidding process converges in finite steps.

For result: 
1) \emph{Incentive}. Prosumers are willing to participate in sharing, and more participants lead to better performance.
2) \emph{Effective}. The equilibrium satisfies physical constraints.
3) \emph{Meaningful}. The price indicates the value of production-demand balance.
4) \emph{Flexible}. Prosumer's role as a seller or buyer is endogenously given instead of predetermined.
5) \emph{Economical}. The total energy sharing cost 
is lower than total self-sufficiency cost and close to the social optimum.

To meet these requirements, we propose an energy sharing mechanism among prosumers. The basic setting of the proposed mechanism is developed in Section III, which is characterized as a generalized Nash game. Main properties of the  generalized Nash equilibrium are revealed. A bidding process is presented in Section IV with proof of its convergence. Simulations results are shown in Section V to validate our findings. Section VI concludes this work.

\section{Energy Sharing Game}
\subsection{Basic Settings}
\begin{figure}[h]
	\centering
	\includegraphics[width=1.0\columnwidth]{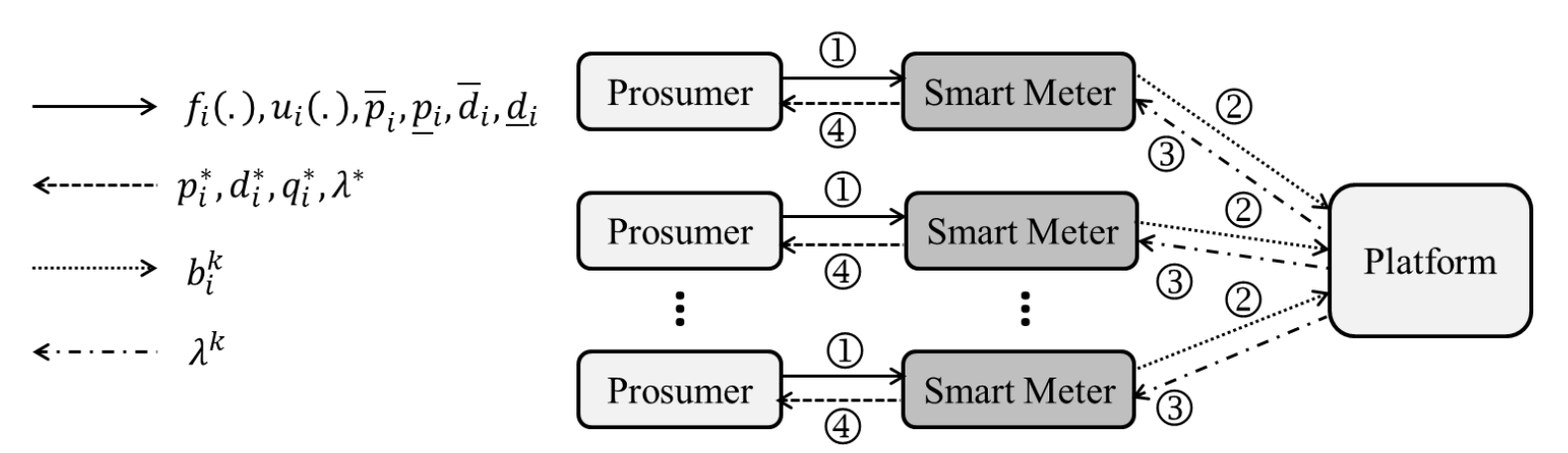}
	\caption{Energy sharing framework between prosumers and the platform.}
	\label{fig:energy-sharing}
	\vspace{-1mm}
\end{figure}

Prosumers participate in an energy sharing market to exchange energy with each other and make individual decisions to maintain energy balancing. 
Specifically, prosumer $i$ imports net energy $q_i$ at market clearing price $\lambda$, which means paying $\lambda q_i$ to buy energy if $q_i>0$ and otherwise receiving revenue $-\lambda q_i$ by selling energy. 
The sharing framework is shown in Fig. \ref{fig:energy-sharing}. Each prosumer is connected to a platform via a smart meter through a bidirectional information channel. The information flow is explained below.

\textbf{Step 1:} (Initialization) Each prosumer $i$ enters its private parameters $f_i(.)$, $u_i(.)$, $\underline{p}_i$, $\overline{p}_i$, $\underline{d}_i$, $\overline{d}_i$ to its smart meter $i$. Set $\lambda^1=0$, and $k=1$. Choose tolerance $\epsilon$.

\textbf{Step 2:} Each smart meter $i$ updates its bid $b_i^{k+1}$ based on the latest $\lambda^k$, and sends it to the platform.

\textbf{Step 3:} After receiving all the bids $b_i^{k+1},\forall i \in \mathcal{I}$, the platform updates price $\lambda^{k+1}$ and sends it back to all the smart meters.

\textbf{Step 4:} If $|\lambda^{k}-\lambda^{k+1}| \leq \epsilon$, $\lambda^*=\lambda^{k+1}$, go to \textbf{Step 5}; otherwise, $k=k+1$ and go to \textbf{Step 2}.

\textbf{Step 5:} Each smart meter determines the optimal production $p_i^*$, demand $d_i^*$, and sharing quantity $q_i^*$ based on $\lambda^*$, and sends them back to the corresponding prosumer to execute. 

In the procedure above, private information is only required by each prosumer’s own smart meter so that its privacy is well preserved.
Details about the mechanism will be explained in Section IV. The key to our mechanism design is to determine price $\lambda$ and quantity $q_i,\forall i \in \mathcal{I}$ based on prosumers' bids $b_i,\forall i \in \mathcal{I}$. We use the generalized demand (or supply) function \cite{hobbs2000strategic} to depict their relationship:
\begin{align}
\label{eq:supply-demand function}
q_i = - a\lambda + b_i, \quad \forall i \in \mathcal{I}
\end{align}
where $a>0$ is a parameter for market sensitivity, and $b_i$ is prosumer $i$'s bid. Market clearing requires $\sum_{i \in \mathcal{I}} q_i=0$ for power balance. Therefore, $ \sum_{i \in \mathcal{I}}(-a\lambda+b_i)=-aI\lambda+\sum_{i \in \mathcal{I}}b_i=0$
and the price turns out to be
\begin{align}
\label{eq:clearing}
\lambda=\frac{\sum_{i \in \mathcal{I}} b_i}{aI}.
\end{align}

Equation \eqref{eq:supply-demand function} is from the typical demand curve where quantity $q_i$ is decreasing with price $\lambda$ \cite{hobbs2000strategic}. We extend it by allowing $q_i$ to be negative for selling energy.
The bid $b_i$ indicates prosumer $i \in \mathcal{I}$’s willingness to buy energy. 
Specifically, equations \eqref{eq:supply-demand function}-\eqref{eq:clearing} imply $q_i = b_i - (\sum_{j\in\mathcal{I}} b_j /I)$, which means prosumer $i$ is a buyer ($q_i>0$) if it is more willing to buy than the average, and a seller otherwise.

The objective of each prosumer $i \in \mathcal{I}$ is to minimize its cost of production minus utility of consumption plus the payment for buying energy (or minus the revenue from selling energy), subject to energy balance and capacity limits for production and consumption. Formally:
\bsq
\label{eq:sharing game}
\begin{align}
\mathop{\min}_{p_i,d_i,b_i} ~ & \Gamma_i(p,d,b)\!:=\!f_i(p_i)\!-\!u_i(d_i)\!+\!\lambda(b)(-a\lambda(b)\!+\!b_i) \label{eq:sharing game.1}\\
\mbox{s.t.}~ & p_i -a\lambda(b)+b_i = d_i \label{eq:sharing game.2}\\
~ & \underline{p}_i \le p_i \le \overline{p}_i \label{eq:sharing game.3}\\
~ & \underline{d}_i \le d_i \le \overline{d}_i \label{eq:sharing game.4} \\
~ & \lambda(b) = \frac{\sum_{j \in \mathcal{I}} b_j}{aI} \label{eq:market clearing}
\end{align}
\esq

The proposed mechanism \eqref{eq:sharing game} can be modeled as a game with the following elements: 1) a set of players $\mathcal{I}$; 2) action sets $S_{i}(p_{-i},d_{-i},b_{-i}),\forall i \in \mathcal{I}$ and strategy space $S=\prod_{i \in \mathcal{I}}S_i$; 3) cost functions $\Gamma_i(p,d,b),\forall i \in \mathcal{I}$. We denote the game compactly as $\mathcal{G}=\{\mathcal{I},S,\Gamma\}$. The action of player $i \in \mathcal{I}$ is composed of production $p_i$, consumption $d_i$, and bid $b_i$, with action set $S_i$ defined by \eqref{eq:sharing game.2}-\eqref{eq:market clearing}. A uniform price $\lambda(b)$ is determined by \eqref{eq:market clearing}, which couples all the players and thus depends on other players’ actions $b_{-i}:=(b_j,\forall j \ne i)$. Since $\lambda(b)$ enters 
constraint \eqref{eq:sharing game.2}, action set $S_i$ defined by \eqref{eq:sharing game.2} (with other constraints) also depends on $b_{-i}$, so that the proposed mechanism is a generalized Nash game, whose equilibrium is harder to analyze than a standard Nash game \cite{facchinei2010generalized}.

\begin{definition} A profile $(\hat{p},\hat{d},\hat{b})\in S$ is a generalized Nash equilibrium (GNE) of the sharing game $\mathcal{G}$, if $\forall i \in \mathcal{I}$: 
	\begin{eqnarray}\nonumber
	(\hat{p}_i\!,\hat{d}_i\!,\hat{b}_i)\in \mbox{argmin}~  \Gamma_i(p_i\!,d_i\!,b_i\!,\hat{p}_{-i}\!,\hat{d}_{-i}\!,\hat{b}_{-i}),
	\mbox{s.t.}~\eqref{eq:sharing game.2}-\eqref{eq:market clearing} 
	\end{eqnarray}
\end{definition}

\subsection{Properties of the Sharing Equilibrium}\label{subsec:GNE-properties}

We next unveil three major properties possessed by the equilibrium of the proposed mechanism. Proposition \ref{prop1} affirms existence of an \emph{effective} market equilibrium that satisfies its defining constraints; Proposition \ref{prop2} states that the equilibrium provides adequate \emph{incentive} for prosumers to participate. Proposition \ref{prop3} shows that the proposed mechanism is \emph{economical}, i.e., the total net cost $\sum_{i \in \mathcal{I}} [f_i(\hat p_i)-u_i(\hat d_i)]$ at equilibrium 
approaches the socially optimal net cost $\sum_{i \in \mathcal{I}} [f_i(\tilde{p}_i)-u_i(\tilde{d}_i)]$.

\begin{proposition}
	\label{prop1}
	(Existence and Partial Uniqueness) A GNE of game $\mathcal{G}$ exists if and only if A1 holds. Moreover, for any GNE $(\hat{p},\hat{d},\hat{b})$, the point $(\hat{p},\hat{d})$ is the \emph{unique} optimal solution to:
	\bsq
	\label{eq:central}
	\begin{align}
	\mathop{\min}_{p_i,d_i,\forall i \in \mathcal{I}}~ & \sum \limits_{i=1}^I f_i(p_i) - \sum \limits_{i=1}^I u_i(d_i) +\frac{\sum_{i=1}^I (d_i-p_i)^2}{2a(I-1)} \label{eq:central.1}\\
	\mbox{s.t.}~ & \sum \limits_{i=1}^I p_i = \sum \limits_{i=1}^I d_i : \zeta \label{eq:central.2}\\
	~ & \underline{p}_i \le p_i \le \overline{p}_i : \delta_i^{\pm} , \forall i \in \mathcal{I} \label{eq:central.3}\\
	~ & \underline{d}_i \le d_i \le \overline{d}_i  : \kappa_i^{\pm} , \forall i \in \mathcal{I} \label{eq:central.4}
	\end{align}
	\esq
\end{proposition}

The proof of Proposition \ref{prop1} is in Appendix \ref{apen-1}.  
Note that constraints \eqref{eq:central.2}--\eqref{eq:central.4} and \eqref{eq:opt3.2}--\eqref{eq:opt3.4} are identical, which implies the migration from centralized operation to the distributed mechanism does not sacrifice feasibility. 
Moreover, energy sharing price $\lambda$ equals ``shadow price'' $\zeta$ in \eqref{eq:central}. Different from the unique Nash equilibrium in \cite{chen2019energy}, the GNE here is partially unique. Specifically, $(\hat p, \hat d)$ is unique, but there can be multiple vectors $\hat b$ leading to the same $(\hat p,\hat d)$. The unique vector $(\hat p,\hat d)$ determines the total net cost $\sum_{i \in \mathcal{I}} [f_i(\hat p_i)-u_i(\hat d_i)]$. Therefore, we use “partially unique” instead of “multiple” to highlight the fact that \emph{the market efficiency is uniquely determined.}
Vector $\hat b$ determines prosumer payment $\hat \lambda(-a\hat \lambda+\hat b_i)$ for $i\in\mathcal{I}$. Note the proposed market has a self-balanced budget: $\sum_{i \in \mathcal{I}}\hat \lambda(-a\hat \lambda+\hat b_i)=\hat \lambda \sum_{i \in \mathcal{I}}(-a\hat \lambda+\hat b_i)=0$, so that $\hat b$ only affects the inner profit allocation among prosumers.

Regardless of difference in $\hat b$, prosumers are always incentivized to participate in sharing, as shown by next proposition.
To prepare for it, we define 
a ``self-sufficiency'' problem:
\bsq
\label{eq:ind}
\begin{align}
\mathop{\min}_{p_i,d_i}~& f_i(p_i)-u_i(d_i) \\
\mbox{s.t.}~& p_i=d_i \label{eq:ind.1}\\
~ & \underline{p}_i \le p_i \le \overline{p}_i  \label{eq:ind.2}\\
~ & \underline{d}_i \le d_i \le \overline{d}_i  \label{eq:ind.3}
\end{align}
\esq
and make the following assumptions:

\noindent \textbf{A2}: Problem \eqref{eq:ind} is feasible.

\noindent \textbf{A3}: $J_i(\check{p}_i,\check{d}_i) < 0,\forall i \in \mathcal{I}$.

Under A2, problem \eqref{eq:ind} for each $i\in \mathcal{I}$ has a unique optimal solution $(\check{p}_i,\check{d}_i)$ due to strict convexity of its objective function. 
Assumption A3 reasonably assumes that under self-sufficiency, each prosumer gets a negative net cost (positive net utility).

\begin{proposition}(Pareto improvement)
	\label{prop2}
	Suppose A2 holds, and $(\hat{p},\hat{d},\hat{b})$ is a GNE of game $\mathcal{G}$. We have
	\begin{align}
	\label{eq:pareto}
	J_i(\check{p}_i,\check{d}_i) \ge \Gamma_i(\hat{p},\hat{d},\hat{b}),\forall i
	\end{align}
	where strictly inequality holds for at least one $i\in\mathcal{I}$ unless $(\check{p},\check{d})= (\hat{p},\hat{d})$.
\end{proposition}

The proof of Proposition \ref{prop2} is in Appendix \ref{apen-2}. It verifies that the proposed mechanism can incentivize prosumers to join since no prosumer is worse off and at least one can benefit. A rare special case is that the self-sufficiency solution coincides with the energy sharing equilibrium (in which case it also coincides with the centralized social optimal for problem \eqref{eq:opt3}).

\textbf{Remark}: One possible application scenario of our model and method is the isolated/standalone microgrids \cite{ali2017determination}, which are designed to be energy self-balanced without a grid connection. In a standalone microgrid, the prosumers can be centrally optimized as in problem \eqref{eq:opt3} (whose solution is the social optimum) or work self-sufficiently as in problem \eqref{eq:ind}. It is worth noting that under self-sufficiency, the prosumers may have to sacrifice their utility in order to match its demand with its generation. However, even under the centralized operation, a prosumer does not always outperform what it would be under the self-sufficiency mode. For instance, as shown later in TABLE \ref{tab:result}, Prosumer 2 and 3 actually have lower utility under the centralized operation. To achieve not only a lower social cost but also lower individual costs than self-sufficiency, we propose the energy sharing mechanism to enable exchanges among prosumers within a microgrid. We prove in Proposition \ref{prop2} that all prosumers have the incentives to join energy sharing since none of them will become worse-off, which is one advantage of the proposed mechanism.

Additionally, our model could be extended to a system with grid connection. For example, if we allow prosumer $i \in \mathcal{I}$ to purchase additional net power $p_i^g$ from the grid at price $\lambda^g$, then since demand-side DERs usually have lower production cost than thermal units in the grid, it is reasonable to assume $\dot f_i(p_i)<\lambda^g$ for all $p_i \in [\underline p_i,~ \overline{p}_i]$ and $i \in \mathcal{I}$. If assumptions A1 and A2 stills hold, at the optimal points of \eqref{eq:opt3} and \eqref{eq:ind} there is $p_i^g=0,\forall i \in \mathcal{I}$. Therefore, the results in this paper can be readily applied. Even without assumptions A1 and A2, where prosumers might buy from the grid, prosumer $i \in \mathcal{I}$ can still obtain at least the same net utility in sharing as that in self-sufficiency by letting $b_i=\sum \nolimits_{j \ne i}\bar b_j/(I-1)$. In other words, Proposition \ref{prop2} can be proved following a similar procedure to that in Appendix \ref{apen-2}.

Although prosumers have incentives to share energy, there is still a gap between the total net cost of energy sharing \eqref{eq:sharing game} and the socially optimal net cost for \eqref{eq:opt3}. 
Our next proposition bounds this gap in terms of price-of-anarchy.
\begin{definition}(Price of Anarchy, PoA \cite{johari2011parameterized})
Consider game $\mathcal{G}=\{\mathcal{I},S,\Gamma\}$. Let $S_{eq} \subseteq S$ be the set of strategies in equilibrium. \emph{Price of Anarchy (PoA)} of game $\mathcal{G}$ is the ratio of the total cost between the worst equilibrium and the social optimal:
\begin{align}
    \mbox{PoA}(\mathcal{G}):=\frac{\mathop{\max}_{s \in S_{eq}} \sum \limits_{i=1}^I \Gamma_i(s)}{\mathop{\min}_{s \in S} \sum \limits_{i=1}^I \Gamma_i(s)}
\end{align}
\end{definition}
PoA measures how the overall efficiency of a game degrades due to strategic behaviors of players. Particularly, a PoA equal to $1$ implies the game achieves social optimal.

\begin{proposition}(Tendency)
	\label{prop3} Suppose A1--A3 hold, and $\underline p_i$, $\overline p_i$, $\underline d_i$, $\overline d_i$, $f_i(.)$, $u_i(.)$ over all $i \in \mathcal{I}$ are uniformly bounded by numbers independent from $I$. Given $I$, let $(\hat{p}(I),\hat{d}(I),\hat{b}(I))$ be a GNE of game $\mathcal{G}$, and $(\tilde{p}(I),\tilde{d}(I))$ be the unique optimal solution of \eqref{eq:opt3}. We have
	\begin{align}
	\label{eq:PoA}
	    \mbox{PoA}(\mathcal{G})=\frac{J\left(\hat{p}(I),\hat{d}(I)\right)}{J\left(\tilde{p}(I),\tilde{d}(I)\right)} \ge 1-\frac{C}{I-1}
	\end{align}
where $C$ is a constant. Moreover, there is
\begin{eqnarray}
	\label{eq:converge2}
	 \lim\limits_{I \to \infty} \left|\hat{p}_i(I)-\tilde{p}_i(I)\right| = \lim\limits_{I \to \infty} \left|\hat{d}_i(I)-\tilde{d}_i(I)\right| = 0, ~\forall i \in \mathcal{I}. 
	\end{eqnarray}
\end{proposition}
The proof of Proposition \ref{prop3} is in Appendix \ref{apen-3}.
Note that PoA is conventionally larger than 1 with a positive cost at social optimal \cite{johari2011parameterized}. In our work, by A3 and Proposition \ref{prop2}, the total net cost is consistently negative across self-sufficiency, energy sharing, and centralized socially optimal mechanisms, which makes PoA less than 1. Proposition 3 shows that both the total net cost and prosumer strategies at a GNE of the proposed mechanism converge to those at the centralized social optimal, with an increasing number of participating prosumers.

\textbf{Remark}: PoA is an important concept measuring inefficiency of a market. A common phenomenon is that fiercer competition leads to a more efficient market, but this is not always true. Here is a counter example: There are $I$ agents in a market. Each agent can bid 0 or 1, and its profit depends on other agents’ bids as shown in TABLE I. 
\begin{table}[h]
	\renewcommand{\arraystretch}{1.3}
	\renewcommand{\tabcolsep}{1em}
	\centering
	\caption{Payoff matrix for each agent}
	\label{tab:data1}
	\begin{tabular}{ccc}
		\hline 
		bid & 0 & 1 \\
		\hline
		all other agents’ bids are 0 & 1 & 0\\
		other cases & 4 & 3\\
		\hline
	\end{tabular}
\end{table}

\noindent Given other agents’ bids, the best strategy of an agent is always to bid 0. Therefore, the market equilibrium is that all the agents bid 0, at which the total profit equals $I$. However, for $I$ large enough, the maximum total profit is $4(I-1)$. Therefore, as we introduce more competition by making $I \rightarrow +\infty$, PoA is actually decreasing (worse) and approaching $1/4$, in which case the added competition does not improve market efficiency.
Therefore, we analyze PoA to ensure that no exception as in the counter example above occurs to the proposed mechanism.
Moreover, analyzing PoA also reveals a nontrivial result \eqref{eq:PoA} that the proposed mechanism approaches social optimal at a rate of $O(1/I)$.

\section{Bidding Process}
This section presents a bidding process and a range of market sensitivity $a$ that guarantees convergence of this process to the GNE characterized in Section \ref{subsec:GNE-properties}. 

\subsection{Bidding Process}
The bidding process is shown structurally in Fig. \ref{fig:energy-sharing} and elaborated in \textbf{Algorithm 1}.
The key to this process is for each prosumer to update its bid without knowing other prosumers' actions. Specifically, at $(k+1)^{\mbox{th}}$ iteration, each prosumer $i \in \mathcal{I}$ utilizes the up-to-date price $\lambda^k$ to estimate (due to the fact that $\lambda(p,d)$ is not known exactly) its optimal solution for problem \eqref{eq:sharing game} which is equivalent to:
\bsq
\label{eq:procedure}
\begin{align}
\mathop{\min}_{p_i,d_i} ~& f_i(p_i)-u_i(d_i)+\lambda(p,d)(d_i-p_i) \label{eq:procedure.1}\\
\mbox{s.t.}~ & \underline{p}_i \le p_i \le \overline{p}_i  \label{eq:procedure.2}\\
~ & \underline{d}_i \le d_i \le \overline{d}_i \label{eq:procedure.3}
\end{align}
\esq
Denote this estimated optimal solution as $(p_i^{k+1},d_i^{k+1})$, and the updated bid of prosumer $i$ is $b_i^{k+1} := d_i^{k+1}-p_i^{k+1}+a\lambda^{k}$. Denote the feasible set of problem \eqref{eq:procedure} as $\mathcal{Y}_i$. 
To estimate $(p_i^{k+1},d_i^{k+1})$, instead of simply replacing the term $\lambda(p,d)(d_i-p_i)$ with $\lambda^k (d_i-p_i)$, prosumer $i$ considers the predicted impact of its decision on price $\lambda(p,d)$, by taking the partial derivative of $\lambda(p,d)(d_i-p_i)$ over $p_i$ (similarly for $d_i$):
\begin{align}
\label{eq:impact}
& \frac{\partial \lambda(p,d)(d_i-p_i)}{\partial p_i}\Big|_{\lambda=\lambda^k} \nonumber\\
=~ &	\left[\frac{\partial \lambda(p,d)}{\partial p_i}(d_i-p_i) -\lambda(p,d)\right]\Big|_{\lambda=\lambda^k} \nonumber\\
= ~ & -\frac{d_i-p_i}{(I-1)a}-\lambda^k
\end{align}
where the last equality is because of
\begin{eqnarray}
\lambda(p,d) &=& \frac{(d_i-p_i)+\sum_{j \ne i}b_j}{(I-1)a} \nonumber
\end{eqnarray}
derived from \eqref{eq:sharing game.2}--\eqref{eq:market clearing}.
We obtain from \eqref{eq:impact} the following objective function as a surrogate for \eqref{eq:procedure.1}:
\begin{align}
\label{eq:objective-eq}
    f_i(p_i)-u_i(d_i)+\frac{(d_i-p_i)^2}{2a(I-1)}+\lambda^k(d_i-p_i)
\end{align}
and thus convert \eqref{eq:procedure} to:
\begin{align}
\label{eq:procedure-eq}
    \mathop{\min}_{p_i,d_i}~ \eqref{eq:objective-eq}, ~\forall (p_i,d_i) \in \mathcal{Y}_i
\end{align}

\begin{algorithm}[t]
	\caption{Energy Sharing Bidding}
	\KwIn{input parameters $f_i(.)$, $u_i(.)$, $\underline{p}_i,\overline{p}_i$, $\underline{d}_i,\overline{d}_i$ into each smart meter $i$, tolerance $\epsilon$.}
	\KwOut{energy sharing results $p^*, d^*, b^*, \lambda^*$.}
	\textbf{Initialization:} $\lambda^1=0$, $k=0$\; 
	\Repeat{$|\lambda^{k+1}-\lambda^k| \le \epsilon$}{
		iteration $k++$
		
		\textbf{prosumer update:}
		
		\For{$i=1; i \le I$}
		{ 
			\begin{align}
			(p_i^{k+1}, d_i^{k+1})~ & \mbox{solves problem \eqref{eq:procedure-eq}} \nonumber\\
			b_i^{k+1}~& := d_i^{k+1}-p_i^{k+1}+a\lambda^k \nonumber
			\end{align}
		}
		
		\textbf{platform update:}
		\begin{align}
		\lambda^{k+1}~& := \frac{\sum_{i=1}^I b_i^{k+1}}{aI} \nonumber
		\end{align}
	}
\end{algorithm}
\subsection{Convergence}
We provide the following condition, under which the proposed bidding process can be proved to converge.

\noindent \textbf{A4}: The market sensitivity $a$ satisfies: $$a \ge \frac{2I-4}{I-1}\mbox{sup}\left\{\frac{1}{\ddot{f}_i(p_i)},-\frac{1}{\ddot{u}_i(d_i)}, ~\forall (p_i,d_i) \in \mathcal{Y}_i,~\forall i\in \mathcal{I} \right\}$$

\begin{proposition}
\label{prop4}
    When A1, A4 hold, \textbf{Algorithm 1} converges to a GNE of the energy sharing game $\mathcal{G}$.
\end{proposition}

For proving convergence of the bidding process, we first give the following lemma with its proof in Appendix \ref{apen-4}.
For conciseness, denote $y_i=[p_i,d_i]^T$, $y=[y_1^T,\dots, y_I^T]^T$, and $\mathcal{Y}=\prod_{i \in \mathcal{I}}\mathcal{Y}_i$. Let $h \in \mathbb{R}^{1 \times 2I}$ be a vector with $h_{2i-1}=1$ and $h_{2i}=-1$ for all $i=1\cdots I$. Define
\begin{align}
    \phi(y):= & \sum \limits_{i=1}^I f_i(p_i)-\sum \limits_{i=1}^I u_i(d_i)+\frac{\sum \limits_{i=1}^I (d_i-p_i)^2}{2a(I-1)}-\frac{(\sum \limits_{i=1}^I d_i-\sum \limits_{i=1}^I p_i)^2}{2aI} \nonumber
\end{align}
and $ L(y,\lambda) := \phi(y)-\lambda hy$ with $\textbf{dom} ~ L = \mathcal{Y}\times \mathbb{R}$.

\begin{lemma}
\label{lemma1}
When A4 holds, $\phi(y)$ is a convex function, and $L(y,\lambda)$ has a (not necessarily unique) saddle point.
\end{lemma}

With Lemma \ref{lemma1}, we next prove Proposition \ref{prop4}.

\begin{proof}
Substituting $b_i^{k+1} := d_i^{k+1}-p_i^{k+1}+a\lambda^{k}$ into $\lambda^{k+1} = (\sum_{i=1}^I b_i^{k+1})/(aI)$, the $k$-th iteration of \textbf{Algorithm 1} becomes: 
\begin{align}
    y_i^{k+1} =~ &  \mbox{argmin} \{\eqref{eq:objective-eq}| y_i \in  \mathcal{Y}_i\},\forall i \in \mathcal{I} \label{eq:iterate.1} \\
    \lambda^{k+1}= ~ & \lambda^k - \frac{hy^{k+1}}{aI}\label{eq:iterate.2}
\end{align} 

Equation \eqref{eq:iterate.1} can be further represented as
\begin{align}
\label{eq:iterate-eq}
    y^{k+1} = ~ & \mbox{argmin}\{\phi(y)-\lambda^khy+\frac{1}{2aI}y^Th^Thy| y \in \mathcal{Y}\}
\end{align}

Utilizing variational inequality and convexity of $\phi(.)$, $y^{k+1} \in \mathcal{Y}$ generated by \eqref{eq:iterate-eq} satisfies
\begin{align}
\label{eq:variational-ineq.1}
   \forall y \in \mathcal{Y}, ~ & \phi(y)-\phi(y^{k+1}) \nonumber\\
  ~&  +(y-y^{k+1})^T\left\{-\lambda^kh^T+\frac{1}{aI}h^T(hy^{k+1})\right\} \ge 0 
\end{align}

Substituting \eqref{eq:iterate.2} into \eqref{eq:variational-ineq.1}, we get
\begin{align}
\label{eq:variational-ineq.2}
    \forall y \in \mathcal{Y}, \phi(y)-\phi(y^{k+1})+(y-y^{k+1})^T(-\lambda^{k+1}h^T) \ge 0
\end{align}

Combining \eqref{eq:variational-ineq.2} and \eqref{eq:iterate.2} gives the following inequality:
\begin{align}
\label{eq:variational-ineq.3}
 &   \left(                
  \begin{array}{c}   
    y-y^{k+1} \\ 
    \lambda-\lambda^{k+1} \\ 
  \end{array}
\right)^T
\left\{
\left(                
  \begin{array}{c}   
    -\lambda^{k+1}h^T \\ 
    hy^{k+1} \\ 
  \end{array}
\right)+
\left(                
  \begin{array}{c}   
    0 \\ 
    aI (\lambda^{k+1}-\lambda^k)\\ 
  \end{array}
\right)
\right\} \nonumber\\
& +    \phi(y)-\phi(y^{k+1})  \ge 0,\forall (y,\lambda) \in \mathcal{Y}\times \mathbb{R}
\end{align}

According to Lemma \ref{lemma1}, let $(y^*,\lambda^*)$ be a saddle point of $L(y,\lambda)$, then we have for any $(y,\lambda) \in \mathcal{Y}\times \mathbb{R}$
\begin{align}
\label{eq:variational-ineq.4}
    \phi(y)-\phi(y^*)+
    \left(                
  \begin{array}{c}   
    y-y^* \\ 
    \lambda-\lambda^*\\ 
  \end{array}
\right)^T
F(y^*,\lambda^*) \ge 0
\end{align}
where the mapping $F(y,\lambda):=[-\lambda h, hy]^T$ and is monotone. 
\footnote{A mapping $F(\lambda,y)$ is monotone if for any $(y_1,\lambda_1), (y_2,\lambda_2) \in \mathcal{Y} \times \mathbb{R}$: 
\begin{eqnarray}\nonumber
\left(\begin{array}{c}   
    y_1-y_2 \\ 
    \lambda_1-\lambda_2\\ 
  \end{array}\right)^T
\left( F(y_1,\lambda_1) - F(y_2, \lambda_2) \right) \geq 0  
\end{eqnarray}}

Since \eqref{eq:variational-ineq.3} holds for all $(y,\lambda)$ in $\mathcal{Y}\times \mathbb{R}$, and particularly for $(y^*,\lambda^*)$, we have
\begin{align}
    ~ & (\lambda^{k+1}-\lambda^*)(\lambda^k-\lambda^{k+1}) \nonumber\\
    \ge~ & \frac{1}{aI} \left\{
     \left(                
  \begin{array}{c}   
    y^{k+1}-y^* \\ 
    \lambda^{k+1}-\lambda^*\\ 
  \end{array}
\right)^T
F(y^{k+1},\lambda^{k+1})+\phi(y^{k+1})-\phi(y^*)
    \right\}
\end{align}
and similarly for \eqref{eq:variational-ineq.4} we have
\begin{align}
    \phi(y^{k+1})-\phi(y^*)+\left(                
  \begin{array}{c}   
    y^{k+1}-y^* \\ 
    \lambda^{k+1}-\lambda^*\\ 
  \end{array}
\right)^TF(y^*,\lambda^*) \ge 0
\end{align}
By monotonicity of mapping $F$, we have
\begin{align}
    (\lambda^{k+1}-\lambda^*)(\lambda^k-\lambda^{k+1}) \ge 0
\end{align}
which implies 
\begin{align}
\label{eq:converge-condition}
    |\lambda^{k+1}-\lambda^*|^2 \le |\lambda^k-\lambda^*|^2 - |\lambda^{k}-\lambda^{k+1}|^2
\end{align}
For every saddle point $( y^*, \lambda^*)$ inequality \eqref{eq:converge-condition} holds. Denote the set of $\lambda^*$ as $\mathcal{W}$. The term $|\lambda^k- \lambda^*|^2$ decreases in each iteration by an amount $|\lambda^{k}-\lambda^{k+1}|^2$, so the sequence $\{|\lambda^k-\lambda^*|^2\}$ converges and the sequence $\{\lambda^k\}$ is bounded. With \eqref{eq:variational-ineq.3} we know that every cluster point of  $\{\lambda^k\}$ belongs to $\mathcal{W}$. With \eqref{eq:converge-condition}, the sequence $\{\lambda^k\}$ only has one cluster point, and thus $\{\lambda^k\}$ converges to a point $ \lambda^* \in \mathcal{W}$.  Substituting $\lambda^*$ into \eqref{eq:iterate.1}, we get $p^k \to p^*$, $d^k \to d^*$, and thus $b^k \to b^*$.

Note that the saddle point $(y^*,\lambda^*)$ of $L(y,\lambda)$ corresponds to a primal-dual optimal of problem \eqref{eq:central}. Since problem \eqref{eq:central} has a unique primal optimal $(\hat p,\hat d)$, we have $p^*=\hat p $, $ d^*=\hat d$; moreover, $\lambda^*=\hat \zeta$ is a dual optimal. Therefore, $(p^*, d^*, b^*)$ is a GNE of the energy sharing game $\mathcal{G}$.
\end{proof}

Proposition \ref{prop4} 
offers a guidance for selecting parameter $a$ to implement the proposed mechanism. 
It also verifies that $\{\lambda^k\}$ converges to the ``shadow price'' $\hat \zeta$ of problem \eqref{eq:central}; moreover, as shown in the proof of Proposition \ref{prop3}, $\hat \zeta$ approaches $\tilde{\lambda}_m$ as $I \to \infty$. Therefore, the energy sharing price is \emph{meaningful} by measuring the value of production-consumption balance.

Assumption A4 is practical. The simulation in Section V-B shows that the bidding process converges with a wide range of $a$, even for some cases where $a$ violates A4 (which is a sufficient but not necessary condition for convergence). 
Although we call $a$ market sensitivity, it is indeed a parameter in the set rule for market clearing. We can adjust $a$ to satisfy A4, in which case the prosumers still have incentives to participate in energy sharing, as claimed by Proposition 2.

\textbf{Remark}: The proposed bidding process falls in the general category of dual gradient method. Therefore, if prosumers do not update their bids at every time step, the process can be modeled in a similar way as a partially asynchronous gradient algorithm, whose convergence can be proved by \cite[Section 7.5]{bertsekas1989parallel}, under certain conditions such as boundedness of time steps during which a prosumer keeps missing its update.

\subsection{Prosumer Rationality and Economic Intuition}
First, by Proposition \ref{prop4}, the bidding process converges to a GNE of game \eqref{eq:sharing game}, at which the market is cleared: $\sum_{i \in \mathcal{I}}q_i = \sum_{i \in \mathcal{I}}(-a\lambda+b_i)=0$, and each prosumer achieves power balance: $p_i+q_i=p_i-a\lambda+b_i=d_i$. 
To generate such a reasonable outcome, we assume the market is executed only after the bidding process converges. 
Second, during the bidding, each prosumer solves \eqref{eq:procedure-eq} whose objective is:
\begin{eqnarray}
\mathop{\min}_{p_i,d_i} ~  f_i(p_i)-u_i(d_i)+\left(\lambda^k+\frac{d_i-p_i}{2a(I-1)}\right)(d_i-p_i) \nonumber
\end{eqnarray}
where $\lambda^k$ is the market announced price for the current iteration and $(d_i-p_i)/(2a(I-1))$ is the predicted impact of prosumer $i$’s decision on price, so that $\lambda^k+(d_i-p_i)/(2a(I-1))$ is prosumer $i$’s predicted price for the next iteration. 
The term $(d_i-p_i)$ following the predicted price is prosumer $i$’s unmet demand which it needs to buy from the market. In summary, prosumer $i$ produces $p_i$, consumes $d_i$, and buys $(d_i-p_i)$ from the market, and its rationality is to minimize its own net cost (production cost - utility + purchase cost from market) while considering its impact on market price. 
Third, 
even though the self power balance constraint \eqref{eq:sharing game.2} is not explicitly in \eqref{eq:procedure-eq}, it is satisfied at equilibrium as each prosumer consistently implements $b_i^{k+1} = d_i^{k+1}-p_i^{k+1}+a\lambda^k$ over iterations.

With the proposed bidding process, we also have an intuitive explanation for Proposition \ref{prop3}. As said, prosumers update their bids considering their impact on price $\lambda$. 
When there is a small number $I$ of prosumers, they constitute a \emph{monopolistic competition} market, where the impact of each prosumer's strategy on price cannot be neglected. When $I$ is large enough, the market is close to \emph{perfectly competitive}, and each prosumer has an infinitesimal influence on price $\lambda$, which can be regarded as exogenously given. In this case, \eqref{eq:impact} reduces to:
\begin{align}\label{eq:impact-reduced}
    \frac{\partial \lambda(d_i-p_i)}{\partial p_i}|_{\lambda=\lambda^k}=-\lambda^k
\end{align}
Following a similar procedure to the proof of Proposition \ref{prop4}, we can show that as $I\rightarrow \infty$, the bidding process with \eqref{eq:impact-reduced} converges to the optimal solution of problem \eqref{eq:opt3}, which is the second statement of Proposition \ref{prop3}.

\section{Simulation}
Numerical experiments are conducted to validate theoretical results. We first run a simple three-prosumer case to verify convergence of the bidding process and efficiency of GNE.

\subsection{Simple Example with Three Prosumers}
\begin{figure}[h]
	\centering
	\includegraphics[width=0.95\columnwidth]{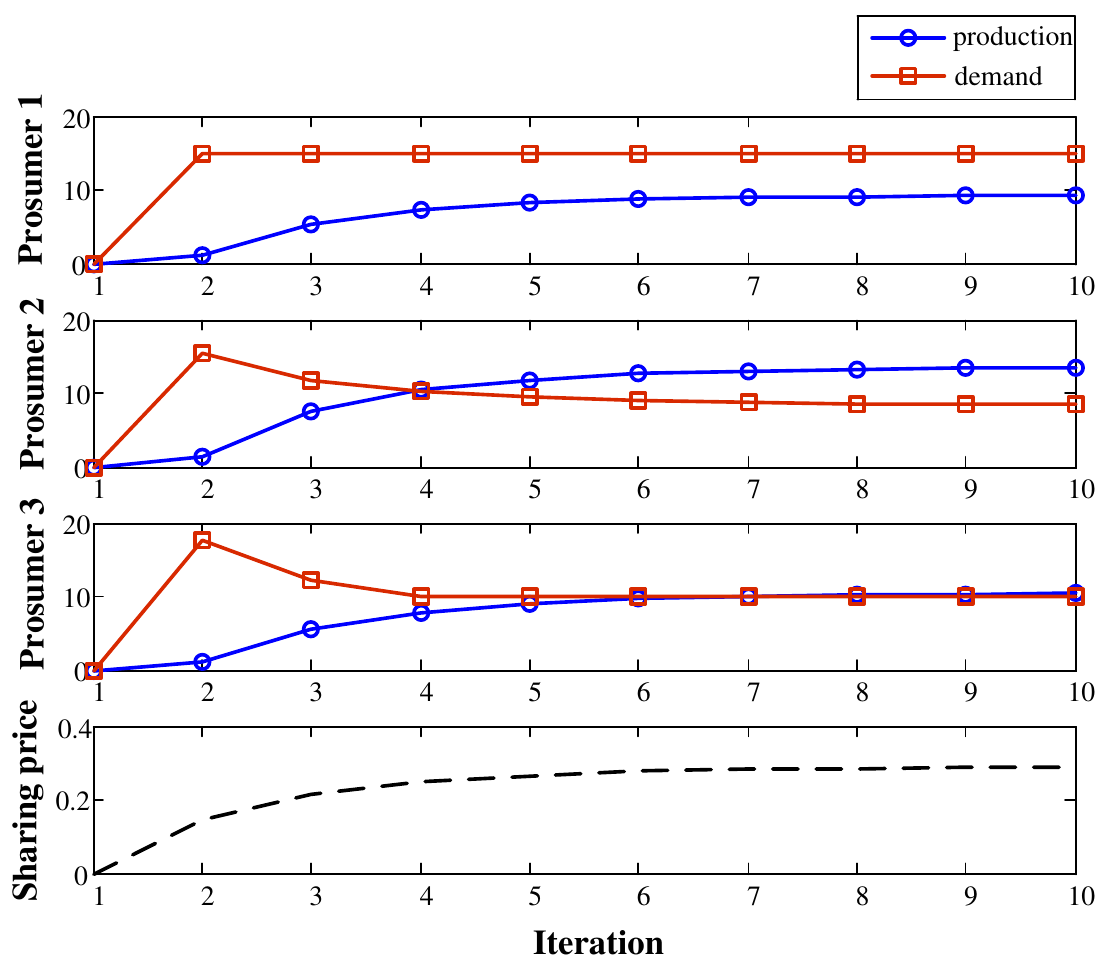}
	\caption{Prosumers' strategies and sharing price over iterations.}
	\label{fig:convergence}
\end{figure}

In the three-prosumer case, market sensitivity is set at $a=100$, cost functions are $f_i(.):=\alpha_i^1p_i^2+\alpha_i^2p_i$, and utility functions are $u_i(.):=\beta_i^1 d_i^2+\beta_i^2 d_i$, where $\alpha_i^1,\alpha_i^2,\beta_i^1,\beta_i^2$, $\forall i \in \{1,2,3\}$ and other parameters are given in TABLE \ref{tab:data1}--\ref{tab:data2}. 
The bidding process in Section IV is used to seek for a GNE. The $p_i^k, d_i^k,\lambda^k$ over iterations are shown in Fig. \ref{fig:convergence}. 
We observe that prosumer strategies and the energy sharing price converge in about 6 iterations. At GNE, the gap between demand and production of a prosumer needs to be bought from the market.
\vspace{-0.5mm}
\begin{table}[h]
	\renewcommand{\arraystretch}{1.3}
	\renewcommand{\tabcolsep}{1em}
	\centering
	\caption{Cost coefficients of prosumers}
	\label{tab:data1}
	\begin{tabular}{ccccc}
		\hline 
		& $\alpha_i^1$ & $\alpha_i^2$ & $\beta_i^1$ & $\beta_i^2$ \\
		Prosumer & ($\$/\mbox{kWh}^2$) & ($\$/\mbox{kWh}$) & ($\$/\mbox{kWh}^2$) & ($\$/\mbox{kWh}$)\\
		\hline
		1 & 0.015 &	0.038 &	-0.008	& 0.8\\
		2 & 0.008 &	0.047 & -0.014 & 0.5\\
		3 & 0.011 & 0.056 & -0.009 & 0.4\\
		\hline
	\end{tabular}
\end{table}
\begin{table}[h]
	\renewcommand{\arraystretch}{1.3}
	\renewcommand{\tabcolsep}{1em}
	\centering
	\caption{Physical limits of prosumers}
	\vspace{-0.5mm}
	\label{tab:data2}
	\begin{tabular}{ccccc}
		\hline 
		Prosumer & $\underline{p}_i(\mbox{kWh})$ & $\overline{p}_i(\mbox{kWh})$ & $\underline{d}_i(\mbox{kWh})$ & $\overline{d}_i(\mbox{kWh})$\\
		\hline
		1 & 0 &	20 & 5 & 15	\\
		2 & 0 & 25 & 7 & 18\\
		3 & 0 & 30 & 10 & 25\\
		\hline
	\end{tabular}
\end{table}

The $(\hat p,\hat d)$ at GNE, the social optimal $(\tilde{p},\tilde{d})$ solved from \eqref{eq:opt3}, and the self-sufficiency strategy $(\check{p},\check{d})$ solved from \eqref{eq:ind} are compared in TABLE \ref{tab:result}. The net costs of all the prosumers are negative, satisfying Assumption A3. Though the centralized social optimal achieves the highest total net utility \$ 10.98, two prosumers become worse-off compared with self-sufficiency: Prosumer 2’s net utility decreases from \$ 2.33 to \$ 0.68, and Prosumer 3 from \$ 1.44 to \$ 1.39. Therefore, Prosumers 2 and 3 may not have the incentive to participate in the centralized operation. Under the proposed energy sharing mechanism, Prosumer 1's net utility increases from \$ 6.25 to \$ 6.90, Prosumer 2 from \$ 2.33 to \$ 2.59, and Prosumer 3 keeps the same. This verifies Proposition \ref{prop2} and shows superior incentive of the proposed mechanism compared to the centralized operation. Moreover, the relative gap between the social optimal and GNE is only (10.98-10.94)/10.98=0.36\%, which verifies efficiency of the energy sharing mechanism.
\begin{table}[h]
	\renewcommand{\arraystretch}{1.3}
	\centering
	\caption{Comparison of three schemes}
	\label{tab:result}
	\begin{tabular}{ccccc}
		\hline 
		Prosumer & $(\hat p,\hat d)$ & $(\tilde{p},\tilde{d})$ & $(\check{p},\check{d})$  \\
		\hline
		1 & (9.3,15.0) & (8.1,15.0) & (15.0,15.0) \\
		Net cost(\$) & -6.90 & -8.91  & -6.25 \\ 
		2 & (13.6,8.4) & (14.6,7.8) & (10.3,10.3) \\
		Net cost(\$) & -2.59 & -0.68 & -2.33  \\
		3 & (10.5,10.0) & (10.2,10.0) & (10.0, 10.0) \\
		Net cost(\$) & -1.44 & -1.39 & -1.44  \\ 
		Total net cost (\$) & -10.94 & -10.98 & -10.03\\ 
		\hline
	\end{tabular}
\end{table}

\begin{figure}[h]
	\centering
	\includegraphics[width=1.0\columnwidth]{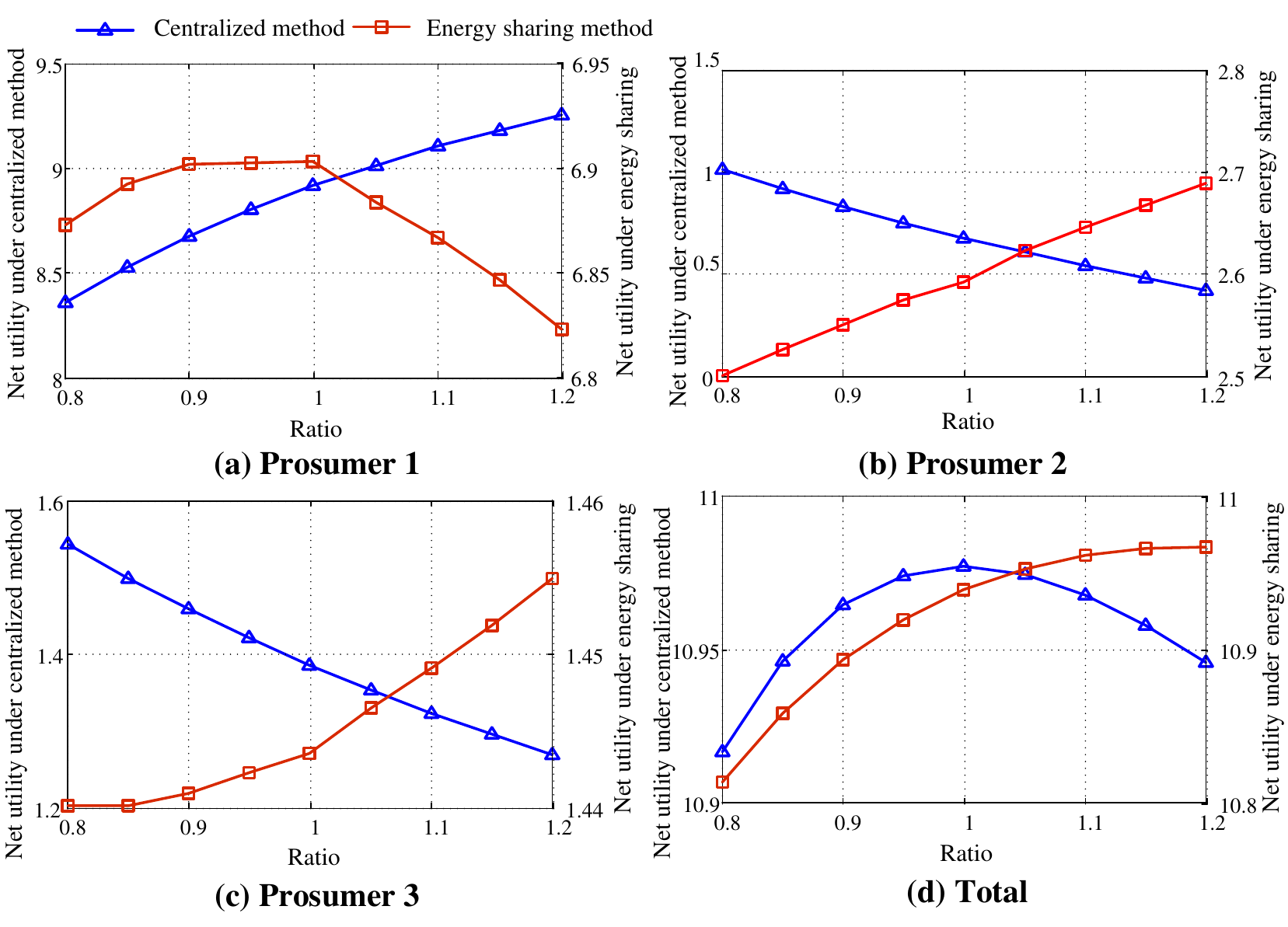}
	\caption{Changes of net utility under the centralized and proposed mechanisms.}
	\label{fig:asymmetric-information}
\end{figure}

We further show the potential of the proposed mechanism in restraining the influence of information asymmetry. Information asymmetry is a crucial problem in market. It describes the situation where a party with more information than others may deliberately misrepresent its information to gain more profit, leading to imbalanced market power or  even market failure \cite{laffont2009theory}. 
We tune the parameter tuple $(\alpha_1^1,\alpha_1^2,\beta_1^1,\beta_1^2)$  from 0.8 to 1.2 times its original value, and show in Fig. \ref{fig:asymmetric-information} the impact of Prosumer 1’s misrepresentation on market equilibrium.

Fig. \ref{fig:asymmetric-information} shows that, under centralized operation, Prosumer 1 tends to report higher  $(\alpha_1^1,\alpha_1^2,\beta_1^1,\beta_1^2)$ to increase its net utility. One consequence, however, is that the net utilities of Prosumers 2 and 3 decline, so does the total net utility of three prosumers. In contrast, under the proposed mechanism, Prosumer 1’s best choice is to report $(\alpha_1^1,\alpha_1^2,\beta_1^1,\beta_1^2)$ truthfully since this leads to its maximum net utility. In this case, information asymmetry does not spoil market equilibrium, which is another merit of the proposed mechanism.

\begin{figure}[h]
	\centering
	\includegraphics[width=0.65\columnwidth]{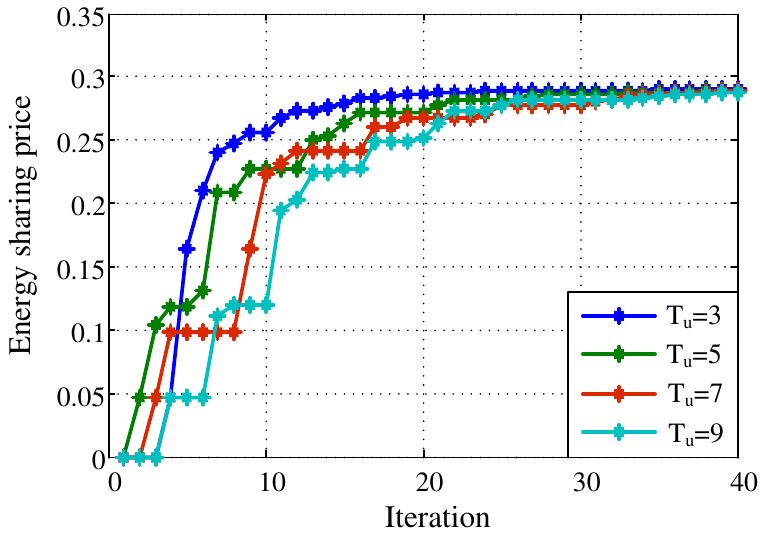}
	\caption{Prices under different upper bounds of update time delay.}
	\label{fig:delay}
\end{figure}

We then illustrate convergence of the bidding process when every prosumer randomly misses its update every iteration with probability 0.8. We change the upper bound of time delay (defined as the number of consecutive iterations during which a prosumer misses its update) from 3 to 9, and the iterates of energy sharing prices are recorded in Fig.\ref{fig:delay}. The energy sharing prices converge under all the different upper bounds of time delay, which indicates that our proposed mechanism is efficient with asynchronous update. Moreover, when a larger upper bound of time delay is allowed, it takes longer time to reach the market equilibrium.

\subsection{Cases with More Prosumers}
\begin{figure}[h]
	\centering
	\includegraphics[width=0.75\columnwidth]{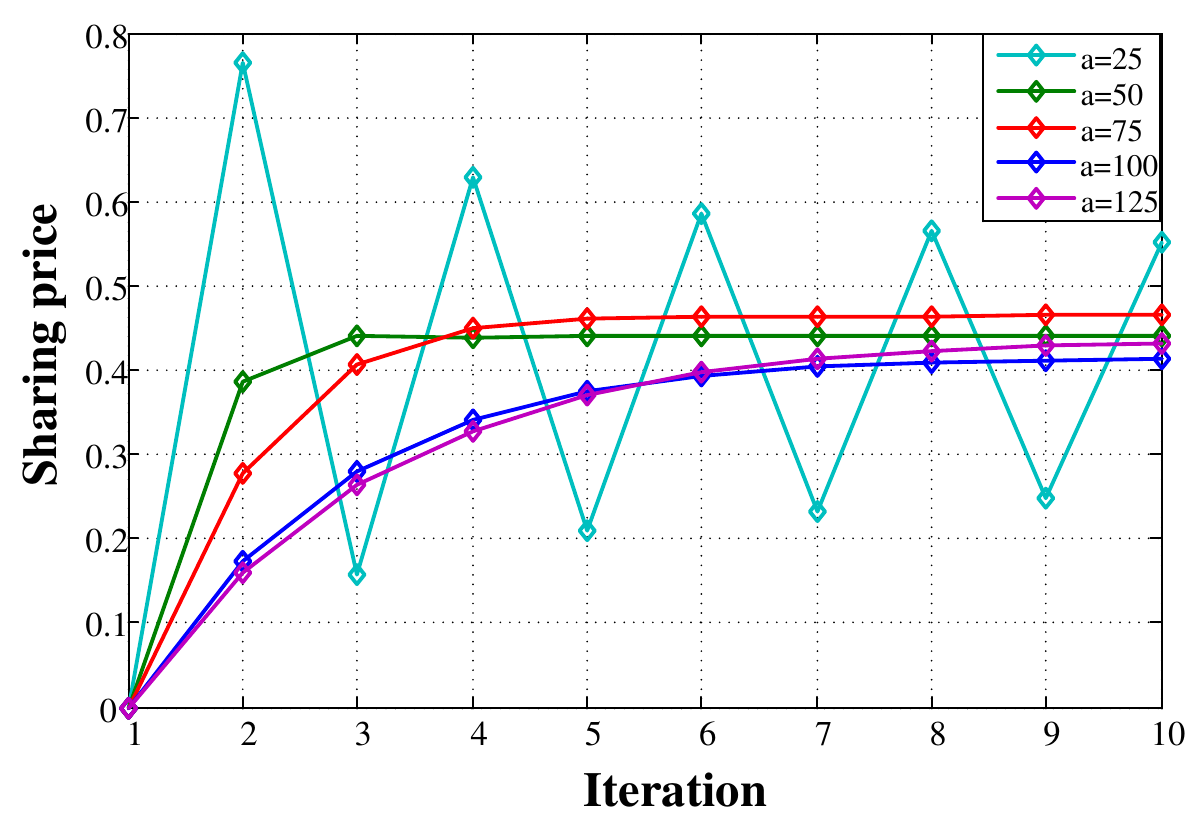}
	\caption{Change of sharing price over iterations, under different $a$.}
	\label{fig:convergenceN}
\end{figure}

We run simulation in a larger case with 50 prosumers to show scalability of the proposed bidding process. Prosumer parameters are uniformly randomly sampled from the following ranges: $\alpha_i^1 \in [0.01,0.02]$, $\alpha_i^2 \in [0.02,0.08]$, $\beta_i^1 \in [-0.01,-0.005]$, $\beta_i^2 \in [0,1]$, $\overline{p}_i \in [20,40]$, $\underline{d}_i \in [5,10]$, $\overline{d}_i \in [15,30]$, and $\underline{p}_i$ is set to zero, $\forall i \in \{1,...,50\}$. We test cases with $a=25$, $50$, $75$, $100$, and $125$. The change of energy sharing price over iterations under each $a$ is plot in a line in Fig. \ref{fig:convergenceN}. When $a=25$, Assumption A4 is violated and the bidding process fails to converge; for other cases, the price converges in about 8 iterations, showing practicability of the proposed process. Note that even for convergent cases $a=50$, $75$, $100$, $125$, A4 is not always met. In other words, A4 is a sufficient but not necessary condition for convergence.

\begin{figure}[h]
	\centering
	\includegraphics[width=0.75\columnwidth]{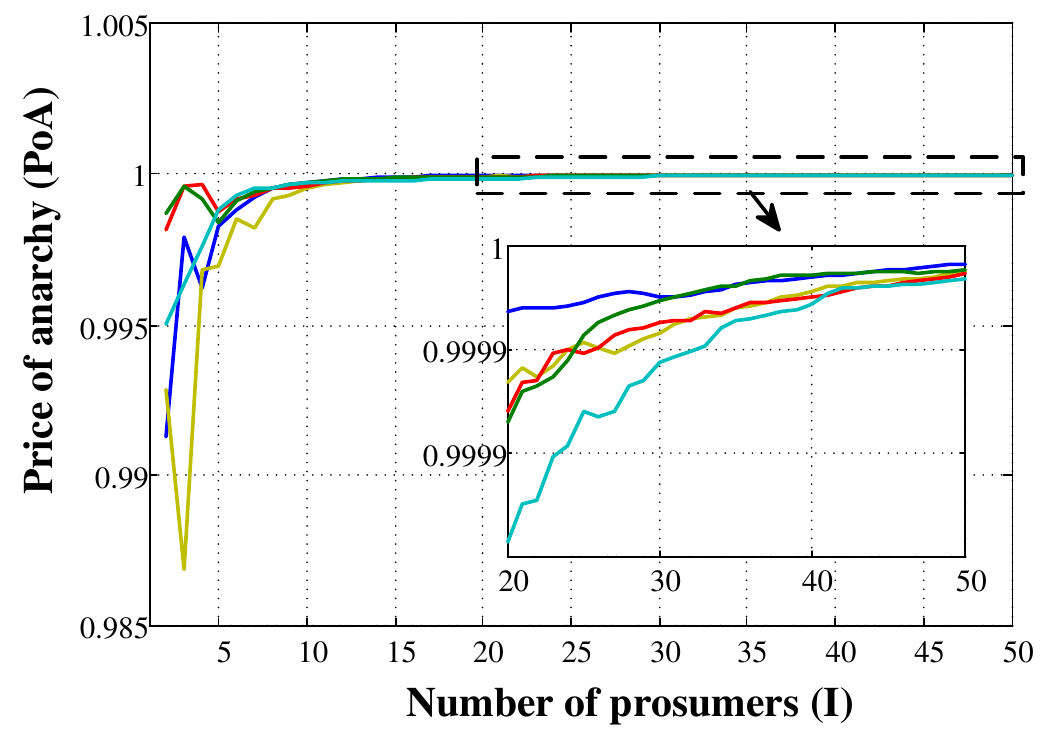}
	\caption{Price of Anarchy (PoA) with an increasing number of prosumers.}
	\label{fig:prosumerN}
\end{figure}

We next test the proposed mechanism with a growing number of prosumers. Increase from 2 to 50 prosumers while selecting parameters in the same way as above and fixing $a=100$. 
The PoA defined in \eqref{eq:PoA} is recorded in Fig. \ref{fig:prosumerN} for five runs (each with a different realization of random parameters). For each case, the PoA converges to 1 as the number of prosumers grows, which validates Proposition \ref{prop3}.

\begin{figure}[h]
	\centering
	\includegraphics[width=0.75\columnwidth]{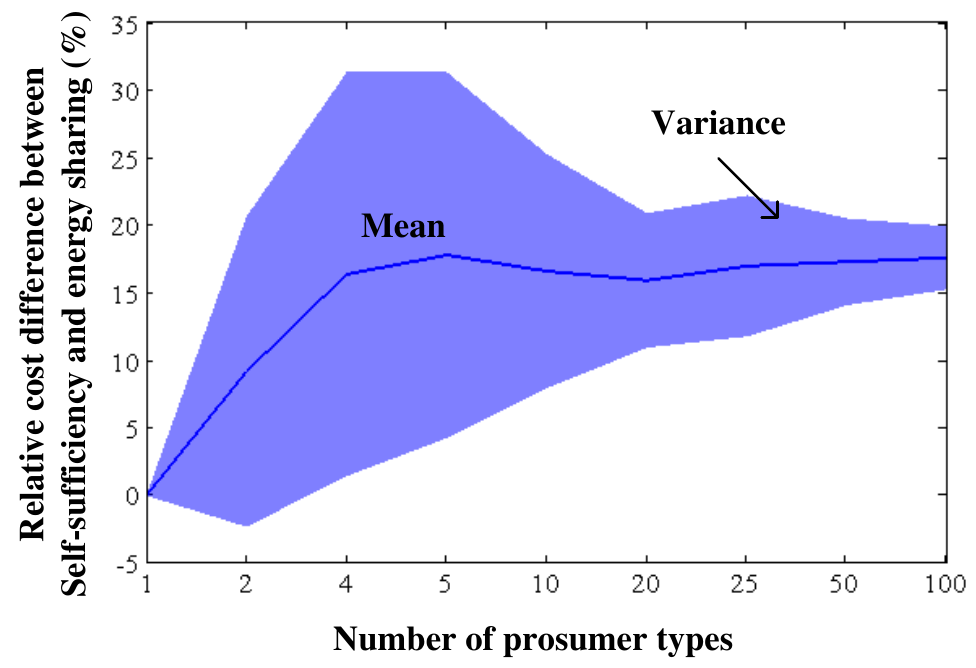}
	\caption{Change of performance of energy sharing with prosumer diversity.}
	\label{fig:prosumer-type}
\end{figure}

We further investigate how prosumer diversity would influence the outcome of energy sharing. The number of prosumers is fixed to 100. At the beginning, all the prosumers have the same parameters, including cost function, utility function, upper/lower bounds $\underline{p}_i,\overline{p}_i,\underline{d}_i,\overline{d}_i,\forall i \in \mathcal{I}$. Then, we gradually add
diversity by increasing the number
of prosumer types \footnote{A prosumer $i \in \mathcal{I}$ in this paper is characterized by four factors: cost function $f_i(.)$, utility function $u_i(.)$, lower and upper bounds of production $\underline{p}_i,\overline{p}_i$, lower and upper bounds of demand $\underline{d}_i,\overline{d}_i$. These four parameters define a “prosumer type”; specifically, prosumers who have the same value of these parameters are called the same type of prosumers.}. Fifty (50) random scenarios are tested for
each degree of diversity, and the mean and variance of the
relative cost difference (saving) of energy sharing versus self-sufficiency
are plotted in Fig.\ref{fig:prosumer-type}. 
With a growing diversity, the mean saving increases, and the variance of saving decays, which demonstrates that more diversified prosumers can lead to more efficient and stable performance of energy sharing.

\section{Conclusion}
We proposed a scalable distributed mechanism for energy sharing to better invoke prosumer flexibility. In the proposed mechanism, a prosumer sends a bid to the market platform without revealing its private information, while its adjustable production and demand and capacity constraints are fully considered. 
The energy sharing mechanism is modeled as a generalized Nash game, whose equilibrium always exists and is partially unique.
At equilibrium, a Pareto improvement is achieved so that every prosumer has the incentive to participate in sharing. By analyzing the price-of-anarchy (PoA), we proved that the performance of energy sharing approaches the centralized social optimal with an increasing number of prosumers. 
A practicable bidding process is presented and its convergence condition is provided. This paper provides insights into market mechanism design in a prosumer era. Future directions include incorporating renewable uncertainties, considering bounded rationality, and characterizing how big data may help improve the performance of energy sharing.

\ifCLASSOPTIONcaptionsoff
\newpage
\fi

\bibliographystyle{IEEEtran}
\bibliography{IEEEabrv,my-bibliography}

\begin{thebibliography}{10}
\providecommand{\url}[1]{#1}
\csname url@samestyle\endcsname
\providecommand{\newblock}{\relax}
\providecommand{\bibinfo}[2]{#2}
\providecommand{\BIBentrySTDinterwordspacing}{\spaceskip=0pt\relax}
\providecommand{\BIBentryALTinterwordstretchfactor}{4}
\providecommand{\BIBentryALTinterwordspacing}{\spaceskip=\fontdimen2\font plus
\BIBentryALTinterwordstretchfactor\fontdimen3\font minus
  \fontdimen4\font\relax}
\providecommand{\BIBforeignlanguage}[2]{{%
\expandafter\ifx\csname l@#1\endcsname\relax
\typeout{** WARNING: IEEEtran.bst: No hyphenation pattern has been}%
\typeout{** loaded for the language `#1'. Using the pattern for}%
\typeout{** the default language instead.}%
\else
\language=\csname l@#1\endcsname
\fi
#2}}
\providecommand{\BIBdecl}{\relax}
\BIBdecl

\bibitem{orrell20162015}
``2017 distributed wind market report,'' Office of {Energy} {Efficiency} \&
  {Renewable} {Energy}, Tech. Rep., 2017.

\bibitem{agnew2015effect}
S.~Agnew and P.~Dargusch, ``Effect of residential solar and storage on
  centralized electricity supply systems,'' \emph{Nature Climate Change},
  vol.~5, no.~4, p. 315, 2015.

\bibitem{parag2016electricity}
Y.~Parag and B.~K. Sovacool, ``Electricity market design for the prosumer
  era,'' \emph{Nature energy}, vol.~1, no.~4, p. 16032, 2016.

\bibitem{liu2017energy}
N.~Liu, X.~Yu, C.~Wang, C.~Li, L.~Ma, and J.~Lei, ``Energy-sharing model with
  price-based demand response for microgrids of peer-to-peer prosumers,''
  \emph{IEEE Transactions on Power Systems}, vol.~32, no.~5, pp. 3569--3583,
  2017.

\bibitem{zhang2013robust}
Y.~Zhang, N.~Gatsis, and G.~B. Giannakis, ``Robust energy management for
  microgrids with high-penetration renewables,'' \emph{IEEE Transactions on
  Sustainable Energy}, vol.~4, no.~4, pp. 944--953, 2013.

\bibitem{liu2017energy2}
N.~Liu, X.~Yu, C.~Wang, and J.~Wang, ``Energy sharing management for microgrids
  with {PV} prosumers: A {Stackelberg} game approach,'' \emph{IEEE Transactions
  on Industrial Informatics}, vol.~13, no.~3, pp. 1088--1098, 2017.

\bibitem{hayes2020co}
B.~P. Hayes, S.~Thakur, and J.~G. Breslin, ``Co-simulation of electricity
  distribution networks and peer to peer energy trading platforms,''
  \emph{International Journal of Electrical Power \& Energy Systems}, vol. 115,
  p. 105419, 2020.

\bibitem{chen2019energy}
Y.~Chen, S.~Mei, F.~Zhou, S.~H. Low, W.~Wei, and F.~Liu, ``An energy sharing
  game with generalized demand bidding: Model and properties,'' \emph{IEEE
  Transactions on Smart Grid}, vol.~11, no.~3, pp. 2055--2066, 2020.

\bibitem{chen2018analyzing}
Y.~Chen, W.~Wei, F.~Liu, Q.~Wu, and S.~Mei, ``Analyzing and validating the
  economic efficiency of managing a cluster of energy hubs in multi-carrier
  energy systems,'' \emph{Applied energy}, vol. 230, pp. 403--416, 2018.

\bibitem{Piclo}
\emph{Piclo website. Available:}, \url{https://piclo.uk/}.

\bibitem{mengelkamp2018designing}
E.~Mengelkamp, J.~G{\"a}rttner, K.~Rock, S.~Kessler, L.~Orsini, and
  C.~Weinhardt, ``Designing microgrid energy markets: A case study: The
  {Brooklyn} microgrid,'' \emph{Applied Energy}, vol. 210, pp. 870--880, 2018.

\bibitem{Enexa}
\emph{{LO3}energy website. Available:},
  \url{https://lo3energy.com/innovations/}.

\bibitem{chakraborty2018sharing}
P.~Chakraborty, E.~Baeyens, K.~Poolla, P.~P. Khargonekar, and P.~Varaiya,
  ``Sharing storage in a smart grid: A coalitional game approach,'' \emph{IEEE
  Transactions on Smart Grid}, vol.~10, no.~4, pp. 4379--4390, 2018.

\bibitem{han2018incentivizing}
L.~Han, T.~Morstyn, and M.~McCulloch, ``Incentivizing prosumer coalitions with
  energy management using cooperative game theory,'' \emph{IEEE Transactions on
  Power Systems}, vol.~34, no.~1, pp. 303--313, 2018.

\bibitem{zhou2020distributed}
Q.~Zhou, M.~Shahidehpour, A.~Paaso, S.~Bahramirad, A.~Alabdulwahab, and
  A.~Abusorrah, ``Distributed control and communication strategies in networked
  microgrids,'' \emph{IEEE Communications Surveys \& Tutorials}, vol.~22,
  no.~4, pp. 2586--2633, 2020.

\bibitem{qi2017sharing}
W.~Qi, B.~Shen, H.~Zhang, and Z.-J.~M. Shen, ``Sharing demand-side energy
  resources-a conceptual design,'' \emph{Energy}, vol. 135, pp. 455--465, 2017.

\bibitem{han2019estimation}
L.~Han, T.~Morstyn, and M.~McCulloch, ``Estimation of the {Shapley} value of a
  peer-to-peer energy sharing game using coalitional stratified random
  sampling,'' \emph{arXiv preprint arXiv:1903.11047}, 2019.

\bibitem{long2019game}
C.~Long, Y.~Zhou, and J.~Wu, ``A game theoretic approach for peer to peer
  energy trading,'' \emph{Energy Procedia}, vol. 159, pp. 454--459, 2019.

\bibitem{liu2019peer}
Y.~Liu, L.~Wu, and J.~Li, ``Peer-to-peer (p2p) electricity trading in
  distribution systems of the future,'' \emph{The Electricity Journal},
  vol.~32, no.~4, pp. 2--6, 2019.

\bibitem{morstyn2018bilateral}
T.~Morstyn, A.~Teytelboym, and M.~D. McCulloch, ``Bilateral contract networks
  for peer-to-peer energy trading,'' \emph{IEEE Transactions on Smart Grid},
  vol.~10, no.~2, pp. 2026--2035, 2018.

\bibitem{ryu2020real}
Y.~Ryu and H.-W. Lee, ``A real-time framework for matching prosumers with
  minimum risk in the cluster of microgrids,'' \emph{IEEE Transactions on Smart
  Grid}, 2020.

\bibitem{ostrovsky2008stability}
M.~Ostrovsky, ``Stability in supply chain networks,'' \emph{American Economic
  Review}, vol.~98, no.~3, pp. 897--923, 2008.

\bibitem{paudel2018peer}
A.~Paudel, K.~Chaudhari, C.~Long, and H.~B. Gooi, ``Peer-to-peer energy trading
  in a prosumer-based community microgrid: A game-theoretic model,'' \emph{IEEE
  Transactions on Industrial Electronics}, vol.~66, no.~8, pp. 6087--6097,
  2018.

\bibitem{dutta2014game}
P.~Dutta and A.~Boulanger, ``Game theoretic approach to offering participation
  incentives for electric vehicle-to-vehicle charge sharing,'' in \emph{2014
  IEEE Transportation Electrification Conference and Expo (ITEC)}.\hskip 1em
  plus 0.5em minus 0.4em\relax IEEE, 2014, pp. 1--5.

\bibitem{anoh2019energy}
K.~Anoh, S.~Maharjan, A.~Ikpehai, Y.~Zhang, and B.~Adebisi, ``Energy
  peer-to-peer trading in virtual microgrids in smart grids: a game-theoretic
  approach,'' \emph{IEEE Transactions on Smart Grid}, vol.~11, no.~2, pp.
  1264--1275, 2019.

\bibitem{tushar2014three}
W.~Tushar, B.~Chai, C.~Yuen, D.~B. Smith, K.~L. Wood, Z.~Yang, and H.~V. Poor,
  ``Three-party energy management with distributed energy resources in smart
  grid,'' \emph{IEEE Transactions on Industrial Electronics}, vol.~62, no.~4,
  pp. 2487--2498, 2014.

\bibitem{zhang2016bidding}
C.~Zhang, J.~Wu, M.~Cheng, Y.~Zhou, and C.~Long, ``A bidding system for
  peer-to-peer energy trading in a grid-connected microgrid,'' \emph{Energy
  Procedia}, vol. 103, pp. 147--152, 2016.

\bibitem{wang2014game}
Y.~Wang, W.~Saad, Z.~Han, H.~V. Poor, and T.~Ba{\c{s}}ar, ``A game-theoretic
  approach to energy trading in the smart grid,'' \emph{IEEE Transactions on
  Smart Grid}, vol.~5, no.~3, pp. 1439--1450, 2014.

\bibitem{le2020peer}
H.~Le~Cadre, P.~Jacquot, C.~Wan, and C.~Alasseur, ``Peer-to-peer electricity
  market analysis: From variational to generalized {Nash} equilibrium,''
  \emph{European Journal of Operational Research}, vol. 282, no.~2, pp.
  753--771, 2020.

\bibitem{facchinei2010generalized}
F.~Facchinei and C.~Kanzow, ``Generalized {Nash} equilibrium problems,''
  \emph{Annals of Operations Research}, vol. 175, no.~1, pp. 177--211, 2010.

\bibitem{shao2016partial}
C.~Shao, X.~Wang, M.~Shahidehpour, X.~Wang, and B.~Wang, ``Partial
  decomposition for distributed electric vehicle charging control considering
  electric power grid congestion,'' \emph{IEEE Transactions on Smart Grid},
  vol.~8, no.~1, pp. 75--83, 2016.

\bibitem{zhang2014congestion}
L.~Zhang, D.~Feng, J.~Lei, C.~Xu, Z.~Yan, S.~Xu, N.~Li, and L.~Jing,
  ``Congestion surplus minimization pricing solutions when lagrange multipliers
  are not unique,'' \emph{IEEE Transactions on Power Systems}, vol.~29, no.~5,
  pp. 2023--2032, 2014.

\bibitem{samadi2012advanced}
P.~Samadi, H.~Mohsenian-Rad, R.~Schober, and V.~W. Wong, ``Advanced demand side
  management for the future smart grid using mechanism design,'' \emph{IEEE
  Transactions on Smart Grid}, vol.~3, no.~3, pp. 1170--1180, 2012.

\bibitem{wei2014robust}
W.~Wei, F.~Liu, S.~Mei, and Y.~Hou, ``Robust energy and reserve dispatch under
  variable renewable generation,'' \emph{IEEE Transactions on Smart Grid},
  vol.~6, no.~1, pp. 369--380, 2014.

\bibitem{li2011optimal}
N.~Li, L.~Chen, and S.~H. Low, ``Optimal demand response based on utility
  maximization in power networks,'' in \emph{2011 IEEE power and energy society
  general meeting}.\hskip 1em plus 0.5em minus 0.4em\relax IEEE, 2011, pp.
  1--8.

\bibitem{li2015storage}
Z.~Li, Q.~Guo, H.~Sun, and J.~Wang, ``Storage-like devices in load leveling:
  Complementarity constraints and a new and exact relaxation method,''
  \emph{Applied Energy}, vol. 151, pp. 13--22, 2015.

\bibitem{wang2016cooperative}
H.~Wang and J.~Huang, ``Cooperative planning of renewable generations for
  interconnected microgrids,'' \emph{IEEE Transactions on Smart Grid}, vol.~7,
  no.~5, pp. 2486--2496, 2016.

\bibitem{choi2016advanced}
J.-Y. Choi, I.-S. Choi, G.-H. Ahn, and D.-J. Won, ``Advanced power sharing
  method to improve the energy efficiency of multiple battery energy storages
  system,'' \emph{IEEE Transactions on Smart Grid}, vol.~9, no.~2, pp.
  1292--1300, 2016.

\bibitem{mondal2015distributed}
A.~Mondal, S.~Misra, and M.~S. Obaidat, ``Distributed home energy management
  system with storage in smart grid using game theory,'' \emph{IEEE Systems
  Journal}, vol.~11, no.~3, pp. 1857--1866, 2015.

\bibitem{chen2020decentralized}
Y.~Chen, T.~Li, C.~Zhao, and W.~Wei, ``Decentralized provision of renewable
  predictions within a virtual power plant,'' \emph{IEEE Transactions on Power
  Systems}, 2020.

\bibitem{pei2016optimal}
W.~Pei, Y.~Du, W.~Deng, K.~Sheng, H.~Xiao, and H.~Qu, ``Optimal bidding
  strategy and intramarket mechanism of microgrid aggregator in real-time
  balancing market,'' \emph{IEEE Transactions on Industrial Informatics},
  vol.~12, no.~2, pp. 587--596, 2016.

\bibitem{hobbs2000strategic}
B.~F. Hobbs, C.~B. Metzler, and J.-S. Pang, ``Strategic gaming analysis for
  electric power systems: An {MPEC} approach,'' \emph{IEEE Transactions on
  Power Systems}, vol.~15, no.~2, pp. 638--645, 2000.

\bibitem{ali2017determination}
L.~Ali and F.~Shahnia, ``Determination of an economically-suitable and
  sustainable standalone power system for an off-grid town in western
  australia,'' \emph{Renewable energy}, vol. 106, pp. 243--254, 2017.

\bibitem{johari2011parameterized}
R.~Johari and J.~N. Tsitsiklis, ``Parameterized supply function bidding:
  Equilibrium and efficiency,'' \emph{Operations Research}, vol.~59, no.~5, pp.
  1079--1089, 2011.

\bibitem{bertsekas1989parallel}
D.~P. Bertsekas and J.~N. Tsitsiklis, \emph{Parallel and distributed
  computation: numerical methods}.\hskip 1em plus 0.5em minus 0.4em\relax
  Prentice hall Englewood Cliffs, NJ, 1989, vol.~23.

\bibitem{laffont2009theory}
J.-J. Laffont and D.~Martimort, \emph{The theory of incentives: the
  principal-agent model}.\hskip 1em plus 0.5em minus 0.4em\relax Princeton
  university press, 2009.

\end{thebibliography}

\appendix
\makeatletter
\@addtoreset{equation}{section}
\@addtoreset{theorem}{section}
\makeatother
\setcounter{equation}{0}  
\renewcommand{\theequation}{A.\arabic{equation}}
\subsection{Proof of Proposition \ref{prop1}}
\label{apen-1}
Given $\bar{b}_{j}, j \ne i$, prosumer $i$'s problem \eqref{eq:sharing game} can be rewritten as \eqref{eq:sharing-eq} by using $p,d$ to represent $\lambda$, $b_i$.
	\bsq
	\label{eq:sharing-eq}
	\begin{align}
	\mathop{\min}_{p_i,d_i} ~~ & f_i(p_i)-u_i(d_i) +\frac{(d_i-p_i)+\sum_{j \ne i}\bar{b}_j}{(I-1)a}(d_i-p_i) \label{eq:sharing-eq.1}\\
	\mbox{s.t.} ~~ & \underline{p}_i \le p_i \le \overline{p}_i : \mu_i^{\pm} \\
	~~ & \underline{d}_i \le d_i \le \overline{d}_i: \eta_i^{\pm}
	\end{align}
	\esq
	and the optimal $b_i$ is given by 
	\begin{align}
	b_i = d_i-p_i + \frac{d_i-p_i+\sum_{j \ne i} \bar{b}_{j}}{I-1}
	\end{align}
	The Hessian matrix of \eqref{eq:sharing-eq.1} is
	\bq
	\begin{bmatrix} \ddot{f}_i(p_i)+\frac{2}{(I-1)a} & -\frac{2}{(I-1)a} \\ -\frac{2}{(I-1)a} & -\ddot{u}_i(d_i)+\frac{2}{(I-1)a} \end{bmatrix} \succ 0 \nonumber
	\eq
	
	So problem \eqref{eq:sharing-eq} is a strictly convex optimization problem, and its KKT condition \eqref{eq:sharing-kkt} is the necessary and sufficient condition  for the optimal solution.
	\bsq
	\label{eq:sharing-kkt}
	\begin{align}
	\dot{f}_i(p_i) - \frac{2(d_i-p_i)+\sum_{j \ne i}\bar{b}_j}{(I-1)a} -\mu_i^{-}+\mu_i^{+} = 0 \label{eq:sharing-kkt.1}\\
	-\dot{u}_i(d_i) +\frac{2(d_i-p_i)+\sum_{j \ne i}\bar{b}_j}{(I-1)a} -\eta_i^{-}+\eta_i^{+} =0 \label{eq:sharing-kkt.2}\\
	0 \le \mu_i^{-} \perp (p_i-\underline{p}_i) \ge 0 \label{eq:sharing-kkt.4}\\
	0 \le \mu_i^{+} \perp (\overline{p}_i-p_i) \ge 0 \label{eq:sharing-kkt.5}\\
	0 \le \eta_i^{-} \perp (d_i-\underline{d}_i) \ge 0 \label{eq:sharing-kkt.6}\\
	0 \le \eta_i^{+} \perp (\overline{d}_i-d_i) \ge 0 \label{eq:sharing-kkt.7}
	\end{align}
	\esq
	
	Then, a profile $(\hat{p},\hat{d},\hat{b})$ is a GNE of $\mathcal{G}$ if and only if $\forall i \in \mathcal{I}$, there exists $\mu_i^{\pm}, \eta_i^{\pm}$, such that $(\hat{p}_i,\hat{d}_i)$ together with $\mu_i^{\pm}, \eta_i^{\pm}$ satisfies \eqref{eq:sharing-kkt} (where $\overline b_j$ is replaced by $\hat b_j$), and $\hat{b}$ satisfies:
	\begin{align}
	\label{eq:condition-b}
	\hat{b}_i = \hat{d}_i - \hat{p}_i + \frac{\hat{d}_i-\hat{p}_i+\sum_{j \ne i}\hat{b}_j}{I-1}, \forall i \in \mathcal{I}
	\end{align}
	
	Problem \eqref{eq:central} is also a strictly convex optimization problem with the KKT condition as in \eqref{eq:central-kkt}.
	\bsq
	\label{eq:central-kkt}
	\begin{align}
	\dot{f}_i(p_i) - \frac{d_i-p_i}{(I-1)a} - \zeta -\delta_i^{-} +\delta_i^{+} = 0 , \forall i \in \mathcal{I} \label{eq:central-kkt:stationarity-f}\\
	-\dot{u}_i(d_i) + \frac{d_i-p_i}{(I-1)a} + \zeta - \kappa_i^{-}+\kappa_i^{+} = 0 , \forall i \in \mathcal{I} \label{eq:central-kkt:stationarity-u}\\
	\sum_{i=1}^I p_i = \sum_{i=1}^I d_i \label{eq:central-kkt:balance}\\
	0 \le \delta_i^{-} \perp (p_i-\underline{p}_i) \ge 0 ,\forall i \in \mathcal{I}\\
	0 \le \delta_i^{+} \perp (\overline{p}_i-p_i) \ge 0,\forall i \in \mathcal{I} \\
	0 \le \kappa_i^{-} \perp (d_i-\underline{d}_i) \ge 0,\forall i \in \mathcal{I} \\
	0 \le \kappa_i^{+}\perp (\overline{d}_i-d_i) \ge 0 ,\forall i \in \mathcal{I}
	\end{align}
	\esq
	
	Suppose a profile $(\hat{p},\hat{d},\hat{b})$ is a GNE of $\mathcal{G}$, and $\hat{\mu}_i^{\pm},\hat{\eta}_i^{\pm},\forall i\in \mathcal{I} $ are the corresponding dual variables, such that \eqref{eq:sharing-kkt}, \eqref{eq:condition-b} are met. Obviously, \eqref{eq:opt3.3} and \eqref{eq:opt3.4} are satisfied. Summing up \eqref{eq:condition-b} for all $i$ gives equation \eqref{eq:opt3.2}. Thus, A1 holds.
	
	Denote
	\bq
	\zeta_i = \frac{\hat{d}_i-\hat{p}_i+\sum_{j \ne i} \hat{b}_j}{(I-1)a},\forall i
	\eq
	
	Condition \eqref{eq:condition-b} indicates that $\forall i \in \mathcal{I}, \hat{d}_i-\hat{p}_i-\hat{b}_i$ are equal, and $\zeta_i$ are equal. Let $\zeta:=\zeta_i,\forall i \in \mathcal{I}$, $\delta_i^{\pm}=\hat{\mu}_i^{\pm}$, $\kappa_i^{\pm}=\hat{\eta}_i^{\pm}$. Then $(\hat{p},\hat{d})$, $\delta^{\pm}$, $\kappa^{\pm}$, and $\zeta$ satisfy the KKT condition \eqref{eq:central-kkt}. Thus, $(\hat{p},\hat{d})$ is the optimal solution of problem \eqref{eq:central} and is unique. 
	
	When A1 holds, problem \eqref{eq:central} is also feasible and has a unique optimal solution $(\hat{p},\hat{d})$ as well as an optimal dual solution ($\hat \delta^{\pm}, \hat \kappa^{\pm}, \hat \zeta$), which together satisfy \eqref{eq:central-kkt}. Let $\mu_i^{\pm}=\hat{\delta}_i^{\pm}$, $\eta_i^{\pm}=\hat{\kappa}_i^{\pm}$, and 
	\begin{align}
	b_i=\hat{d}_i-\hat{p}_i+a\hat{\zeta}
	\end{align}
	Then $(\hat p,\hat d,b)$ and ($\mu^{\pm}, \eta^{\pm}$) satisfy \eqref{eq:sharing-kkt}-\eqref{eq:condition-b}, which implies $(\hat p,\hat d,b)$ is a GNE.

\setcounter{equation}{0}  
\renewcommand{\theequation}{B.\arabic{equation}}
\subsection{Proof of Proposition \ref{prop2}}
\label{apen-2}
	Note that A2 implies A1. For prosumer $i$, given other prosumers' strategies ($\bar{p}_{-i},\bar{d}_{-i},\bar{b}_{-i}$), it can choose $p_i=\check{p}_i$, $d_i=\check{d}_i$ and $b_i=\sum_{j \ne i}\bar{b}_j/(I-1)$, so that $-a\lambda+b_i=0$ and $\Gamma_i(p_i,d_i,b_i,\bar{p}_{-i},\bar{d}_{-i},\bar{b}_{-i})=J_i(\check{p},\check{d})$. Since prosumer $i$ aims at minimizing its net cost at GNE, we have
	\begin{align}\nonumber
	J_i(\check{p}_i,\check{d}_i) \ge \Gamma_i(\hat{p},\hat{d},\hat{b})
	\end{align}
	Suppose \eqref{eq:pareto} holds with equality for all $i$. Adding \eqref{eq:sharing game.1} over all $i\in \mathcal{I}$ leads to:  
	\begin{align}
	\sum \limits_{i=1}^I \Gamma_i(\hat{p},\hat{d},\hat{b}) = \sum \limits_{i=1}^I J_i(\hat{p}_i,\hat{d}_i)
	\end{align}
	Thus, $\sum_{i \in \mathcal{I}} J_i(\check{p}_i,\check{d}_i)=\sum_{i \in \mathcal{I}} J_i(\hat{p}_i,\hat{d}_i) $. The uniqueness of optimal solution of \eqref{eq:central} implies $(\check{p},\check{d})=(\hat{p},\hat{d})$.

\setcounter{equation}{0}  
\renewcommand{\theequation}{C.\arabic{equation}}
\subsection{Proof of Proposition \ref{prop3}}
\label{apen-3}

\noindent\textbf{Part I}: Prove \eqref{eq:PoA}, i.e.,
\begin{eqnarray}\nonumber
	    \mbox{PoA}(\mathcal{G})=\frac{J\left(\hat{p}(I),\hat{d}(I)\right)}{J\left(\tilde{p}(I),\tilde{d}(I)\right)} \ge 1-\frac{C}{I-1}.
\end{eqnarray}

	For simplicity, without causing ambiguity, the $I$ in $(\hat{p}(I),\hat{d}(I))$ and $(\tilde{p}(I),\tilde{d}(I))$ are omitted here. According to Proposition \ref{prop1}, $(\hat{p},\hat{d})$ is the optimal solution of \eqref{eq:central}. Denote $\Omega(p,d):=\sum_{i \in \mathcal{I}} (d_i-p_i)^2$, $\mathcal{S}:=\{(p_i,d_i),\forall i \in \mathcal{I}: \;\mbox{s.t.}\; \eqref{eq:central.2}-\eqref{eq:central.4}\; \mbox{are satisfied.} \}$ Note that $\mathcal S$ is also the feasible set for problem \eqref{eq:opt3}.
	
	For every strategy combination $s = (p,d,b) \in S$, there is:
	\begin{align}
	\sum \limits_{i=1}^I \Gamma_i(p,d,b) = \sum \limits_{i=1}^I J_i(p_i,d_i) = J(p,d)
	\end{align}
	which in particular holds for every GNE $(\hat p, \hat d, \hat b)$ and every $(p^*,d^*,b^*) \in \textnormal{argmin}_{s\in S}\sum_{i=1}^{I} \Gamma_i(s)$. Moreover, one can establish equivalence between the set of all subvectors $(p,d)$ in strategy space $S$ and the feasible set $\mathcal{S}$ of problem \eqref{eq:opt3}, so there must be $(p^*,d^*) = (\tilde p, \tilde d)$.
	Then PoA can be equivalently written as
	\begin{align}
	    \mbox{PoA}(\mathcal{G})=\frac{J(\hat{p},\hat{d})}{J(\tilde{p},\tilde{d})} 
	\end{align}
	
	Obviously $|\Omega(p,d)| \le  C_1 I, \forall (p,d) \in \mathcal{S}$, where
	$$
	C_1:=   \mbox{sup}\left\{|\underline{p}_i-\overline{d}_i|^2,|\overline{p}_i-\underline{d}_i|^2,\forall i \in \mathcal{I}\right\}
	$$
	is independent from $I$ by the uniform bound assumption on $\overline p_i$, $\underline p_i$, $\overline d_i$, $\underline d_i$ for all $i\in\mathcal{I}$. 
	By definition, we have 
	\begin{align}
	J(\tilde{p},\tilde{d}) \le J(\hat{p},\hat{d})
	\end{align}
	and 
	\begin{align}
	J(\tilde{p},\tilde{d})+ \frac{\Omega(\tilde{p},\tilde{d})}{2a(I-1)} \ge J(\hat{p},\hat{d}) + \frac{\Omega(\hat{p},\hat{d})}{2a(I-1)}
	\end{align}
	so that
	\begin{align}
	J(\hat{p},\hat{d}) & \le J(\tilde{p},\tilde{d}) + \frac{\Omega(\tilde{p},\tilde{d})}{2a(I-1)} - \frac{\Omega(\hat{p},\hat{d})}{2a(I-1)} \nonumber\\
	& \le J(\tilde{p},\tilde{d}) + \frac{C_1 I}{a(I-1)}
	\end{align}
	
	When A2 and A3 hold, we have $J(\tilde{p},\tilde{d}) \le J(\check{p},\check{d}) \le C_2 I <0 $, where
	\begin{align}
	    C_2:=  \mbox{sup}\left\{f_i(\check{p}_i)-u_i(\check{d}_i),\forall i \in \mathcal{I}\right\}
	\end{align}
	is independent from $I$ by the uniform bound assumption on $f_i(.)$, $u_i(.)$ for all $i\in\mathcal{I}$. 
	Thus
	\begin{align}
	    1-\mbox{PoA}(\mathcal{G}) = ~ &   \frac{J(\hat{p},\hat{d})-J(\tilde{p},\tilde{d})}{|J(\tilde{p},\tilde{d})|} \nonumber\\
	    \le~ &  \frac{C_1I/a(I-1)}{|C_2| I} \nonumber\\
	    = ~ &  \frac{C}{I-1}
	\end{align}
	where $C:=C_1/(a|C_2|)$.
\\ \newline
 
\noindent\textbf{Part II}: Prove \eqref{eq:converge2}, i.e., 
\begin{eqnarray}\nonumber
	 &&\lim\limits_{I \to \infty} \left|\hat{p}_i(I)-\tilde{p}_i(I)\right| =  \lim\limits_{I \to \infty} \left|\hat{d}_i(I)-\tilde{d}_i(I)\right| = 0,~ \forall i \in \mathcal{I}.
\end{eqnarray}

\noindent\textbf{Sketch of proof.} We notice that the difference between KKT conditions of problems \eqref{eq:opt3} and \eqref{eq:central} only lies in the term $\frac{d_i-p_i}{a(I-1)}$ in stationarity equations \eqref{eq:central-kkt:stationarity-f}--\eqref{eq:central-kkt:stationarity-u}. Due to the uniform bound assumption we made on    
$\underline p_i$, $\overline p_i$, $\underline d_i$, $\overline d_i$, this difference will diminish as prosumer number $I$ increases to infinity. Based on this observation, we can bound the difference between solutions of the two sets of KKT conditions, i.e., between the optimal solutions of problems \eqref{eq:opt3} and \eqref{eq:central}, and show that this difference also diminishes as $I$ increases to infinity. Please see below for a detailed proof. 
\\ \newline
\noindent\textbf{Full proof.} The centralized social optimal problem \eqref{eq:opt3} can be equivalently solved by its KKT condition:\footnote{For convenience, we slightly abuse the notation by denoting capacity-associated dual variables as $(\delta^{\pm}, \kappa^{\pm})$ for both problems \eqref{eq:opt3} and \eqref{eq:central}.} 
   \bsq
	\label{eq:socialoptimal-kkt}
	\begin{align}
	\dot{f}_i(p_i) - \lambda_m -\delta_i^{-} +\delta_i^{+} = 0 , \forall i \in \mathcal{I}\\
	-\dot{u}_i(d_i) +\lambda_m - \kappa_i^{-}+\kappa_i^{+} = 0 , \forall i \in \mathcal{I}\\
	\sum_{i=1}^I p_i = \sum_{i=1}^I d_i \label{eq:socialoptimal-kkt:balance}\\
	0 \le \delta_i^{-} \perp (p_i-\underline{p}_i) \ge 0 ,\forall i \in \mathcal{I}\\
	0 \le \delta_i^{+} \perp (\overline{p}_i-p_i) \ge 0,\forall i \in \mathcal{I} \\
	0 \le \kappa_i^{-} \perp (d_i-\underline{d}_i) \ge 0,\forall i \in \mathcal{I} \\
	0 \le \kappa_i^{+}\perp (\overline{d}_i-d_i) \ge 0 ,\forall i \in \mathcal{I}
	\end{align}
	\esq
All the equations in \eqref{eq:socialoptimal-kkt} except \eqref{eq:socialoptimal-kkt:balance} define the optimal production $p_i$ and consumption $d_i$ in response to a given dual variable $\lambda_m$ as the following functions, for all $i \in \mathcal{I}$:
\begin{eqnarray}\nonumber
&&p_i = \tilde f_i(\lambda_m) :=
\begin{cases}
(\dot{f}_i)^{-1} (\lambda_m),~\textnormal{if}~\dot{f}_i(\underline p_i) < \lambda_m <\dot{f}_i(\overline p_i) \\
\underline p_i, \qquad\qquad  \textnormal{if}~ \lambda_m \leq \dot{f}_i(\underline p_i)  \\
\overline p_i, \qquad \qquad \textnormal{if}~ \lambda_m \geq \dot{f}_i(\overline p_i) 
\end{cases}
\end{eqnarray}    
\begin{eqnarray}\nonumber
&&d_i = \tilde u_i(\lambda_m) :=
\begin{cases}
(\dot{u}_i)^{-1} (\lambda_m),~\textnormal{if}~\dot{u}_i(\overline d_i) < \lambda_m <\dot{u}_i(\underline d_i) \\
\overline d_i, \qquad\qquad  \textnormal{if}~ \lambda_m \leq \dot{u}_i(\overline d_i)  \\
\underline d_i, \qquad\qquad  \textnormal{if}~ \lambda_m \geq \dot{u}_i(\underline d_i) 
\end{cases}
\end{eqnarray}
By our assumptions on $f_i(.)$, $u_i(.)$, for all $i\in \mathcal{I}$, functions $\tilde f_i(.)$ and $-\tilde u_i(.)$ are well defined and monotonically increasing on $\lambda_m \in \mathbb{R}$. By \eqref{eq:socialoptimal-kkt:balance}, $\tilde \lambda_m$ is a dual optimal solution of problem \eqref{eq:opt3} if and only if it solves the following equation:
\begin{eqnarray}\nonumber
&&\sum_{i=1}^{I} \left(\tilde f_i(\tilde \lambda_m) -\tilde u_i (\tilde \lambda_m)\right) =0. 
\end{eqnarray}   

We next look at KKT condition \eqref{eq:central-kkt} which equivalently characterizes primal-dual optimal solutions of problem \eqref{eq:central}. Specifically, all the equations in \eqref{eq:central-kkt} except \eqref{eq:central-kkt:balance} define the optimal $p_i$ and $d_i$ in response to a given dual variable $\zeta$ as functions $f_i^o(\zeta)$ and $u_i^o(\zeta)$, respectively, for all $i \in \mathcal{I}$. Although closed-form expressions of $f_i^o(.)$ and $u_i^o(.)$ are hard to derive, we can establish their relationships with $\tilde f_i(.)$ and $\tilde u_i(.)$, for all $i \in \mathcal{I}$:
\begin{eqnarray}
&& p_i = f_i^o(\zeta) = \tilde f_i \left(\zeta + \frac{d_i - p_i}{a(I-1)}\right) \nonumber \\
&& d_i = u_i^o(\zeta) = \tilde u_i \left(\zeta + \frac{d_i - p_i}{a(I-1)}\right) \nonumber
\end{eqnarray} 
Besides, when $f_i^o(\zeta) \in (\underline p_i, \overline p_i)$ and $u_i^o(\zeta) \in (\underline d_i, \overline d_i)$ are both satisfied, the following equation holds for all $i\in \mathcal{I}$:
\begin{eqnarray}
&&\dot{f}_i\left(f_i^o(\zeta) \right) = \zeta + \frac{u_i^o(\zeta) - f_i^o(\zeta)}{a(I-1)} = \dot{u}_i\left(u_i^o(\zeta) \right) \nonumber
\end{eqnarray}
Taking its derivative over $\zeta$, and combining the cases where capacity constraints are binding, we get for all $i\in\mathcal{I}$:
\begin{eqnarray}\nonumber
&&\dot{f}_i^o(\zeta) = 
\begin{cases}
\frac{1}{\ddot{f}_i(p_i) \left[1+\frac{1}{a(I-1)}\cdot \left( \frac{1}{\ddot{f}_i(p_i)}-\frac{1}{\ddot{u}_i(d_i)} \right) \right] },~\textnormal{if}~\underline p_i < p_i <\overline p_i \\
0, \qquad\qquad  \textnormal{otherwise}  
\end{cases}
\end{eqnarray}    
\begin{eqnarray}\nonumber
&&\dot{u}_i^o(\zeta) = 
\begin{cases}
\frac{1}{\ddot{u}_i(d_i) \left[1+\frac{1}{a(I-1)}\cdot \left( \frac{1}{\ddot{f}_i(p_i)}-\frac{1}{\ddot{u}_i(d_i)} \right) \right] },~\textnormal{if}~\underline d_i < d_i <\overline d_i \\
0, \qquad\qquad  \textnormal{otherwise}  
\end{cases}
\end{eqnarray}
where $p_i=f_i^o(\zeta)$ and $d_i = u_i^o (\zeta)$. By our assumptions on $f_i(.)$, $u_i(.)$, for all $i\in \mathcal{I}$, functions $f_i^o(.)$ and $-u_i^o(.)$ are well defined and monotonically increasing on $\zeta \in \mathbb{R}$. By the power balance constraint \eqref{eq:central-kkt:balance}, $\hat \zeta$ is a dual optimal solution of problem \eqref{eq:central} if and only if it solves the following equation:
\begin{eqnarray}\nonumber
&&\sum_{i=1}^{I} \left(f_i^o(\hat \zeta) -u_i^o (\hat \zeta)\right) =0. 
\end{eqnarray}   

Let $\tilde \lambda_m$ be \emph{any} dual optimal solution of problem \eqref{eq:opt3}, so that $\tilde p_i = \tilde f_i(\tilde\lambda_m)$, $\tilde d_i = \tilde u_i(\tilde \lambda_m)$ for all $i\in \mathcal{I}$ constitute the (unique) primal optimal solution of problem \eqref{eq:opt3}.  
We next show that there must be a dual optimal solution $\hat \zeta$ of problem \eqref{eq:central} which lies near $\tilde \lambda_m$. For that purpose, we denote 
\begin{eqnarray}\nonumber
&& \overline p := \sup\left\{\overline p_i, ~\forall i\in \mathcal{I}\right\},\quad \underline p := \inf\left\{\underline p_i, ~\forall i\in \mathcal{I}\right\} \\
\nonumber
&& \overline d := \sup\left\{\overline d_i, ~\forall i\in \mathcal{I}\right\},\quad \underline d := \inf\left\{\underline d_i, ~\forall i\in \mathcal{I}\right\}
\end{eqnarray}
which all exist and are independent from $I$ by our uniform bound assumption. Define two numbers:
\begin{eqnarray}
&& \zeta^+ := \tilde \lambda_m +\frac{\overline p - \underline d}{a(I-1)}~\geq~
\zeta^- := \tilde \lambda_m +\frac{\underline p - \overline d}{a(I-1)}. \nonumber
\end{eqnarray}

Indeed, there must be
\begin{eqnarray}\nonumber
&&f_i^o(\zeta^+) - u_i^o(\zeta^+) \geq  \tilde f_i(\tilde\lambda_m) - \tilde u_i(\tilde\lambda_m), ~\forall i\in\mathcal{I}
\end{eqnarray} 
which can be verified by assuming $f_i^o(\zeta^+) - u_i^o(\zeta^+) <  \tilde f_i(\tilde\lambda_m) - \tilde u_i(\tilde\lambda_m)$ and deducing a contradiction:
\begin{eqnarray}
&& f_i^o(\zeta^+)= \tilde f_i \left(\zeta^+ +  \frac{d_i^+ - p_i^+}{a(I-1)}\right) 
\geq \tilde f_i \left(\zeta^+ +  \frac{\tilde d_i -\tilde p_i}{a(I-1)}\right) \nonumber 
\\
&& \geq  \tilde f_i \left(\zeta^+ +  \frac{\underline d -\overline p}{a(I-1)}\right) = \tilde f_i(\tilde \lambda_m), \quad \forall i \in \mathcal{I} \nonumber
\end{eqnarray}  
where $p_i^+ = f_i^o(\zeta^+)$, $d_i^+ = u_i^o(\zeta^+)$, for all $i \in \mathcal{I}$. Both inequalities above stem from monotonicity of $\tilde f_i(.)$. Similarly, $-u_i^o(\zeta^+) \geq -\tilde u_i(\tilde \lambda_m)$ for all $i \in \mathcal{I}$, and therefore $f_i^o(\zeta^+) - u_i^o(\zeta^+) \geq  \tilde f_i(\tilde\lambda_m) - \tilde u_i(\tilde\lambda_m)$, which contradicts our assumption. 

We hence further have
\begin{eqnarray}\nonumber
&&\sum_{i=1}^I   \left(f_i^o(\zeta^+) -u_i^o (\zeta^+)\right) \geq  \sum_{i=1}^I   \left(\tilde f_i(\tilde\lambda_m) - \tilde u_i(\tilde\lambda_m)\right) = 0.
\end{eqnarray}
Following the same procedure, we can also show
\begin{eqnarray}
&& \sum_{i=1}^I   \left(f_i^o(\zeta^-) -u_i^o (\zeta^-)\right) \leq   0. \nonumber
\end{eqnarray}
Due to monotonicity of function $\sum_{i=1}^I \left( f_i^o(.) -u_i^o (.) \right)$, there must be $\hat \zeta \in [\zeta^-, \zeta^+]$, such that $ \sum_{i=1}^I   \left(f_i^o(\hat \zeta) -u_i^o (\hat \zeta)\right) =  0$, i.e., $\hat \zeta$ is a dual optimal solution of problem \eqref{eq:central}. Further, $\hat p_i = f_i^o(\hat \zeta)$, $\hat d_i = u_i^o(\hat \zeta)$ for all $i \in \mathcal{I}$ constitute the (unique) primal optimal solution of problem \eqref{eq:central}, which is also the production and consumption profile at GNE.  

To prepare for the final step of our proof, we point out Lipschitz continuity of functions $\tilde f_i(.)$, $\tilde u_i(.)$ for all $i \in \mathcal{I}$. Specifically, due to the uniform bound assumption we made on $\ddot{f}_i(.)$, $\ddot{u}_i(.)$, for all $i \in \mathcal{I}$, there exists a positive constant $\gamma$, which is independent from $I$, such that
\begin{eqnarray}
&& \left|\tilde f_i(x) - \tilde f_i(y)\right| \leq \gamma \left|x-y\right|,~\forall i\in\mathcal{I}, ~\forall x,y \in \mathbb{R}\nonumber\\
  && \left|\tilde u_i(x) - \tilde u_i(y)\right| \leq \gamma \left|x-y\right|,~\forall i\in\mathcal{I}, ~\forall x,y \in \mathbb{R}.\nonumber
\end{eqnarray} 
Denote $\sigma:=\max\left\{|\overline p - \underline d|, ~|\underline p - \overline d|\right\}/a$. For any prosumer number $I$, for all $i \in \mathcal{I}$, we have:
\begin{eqnarray}
\left|\hat p_i - \tilde p_i \right| &\leq& \left|\hat p_i - \tilde f_i(\hat\zeta)\right| + \left|\tilde f_i(\hat \zeta) - \tilde p_i \right| \nonumber\\
&=& \left|\tilde f_i\left(\hat\zeta + \frac{\hat d_i - \hat p_i }{a(I-1)}\right) - \tilde f_i(\hat\zeta)\right| + \left|\tilde f_i(\hat \zeta) - \tilde f_i(\tilde\lambda_m) \right|  \nonumber\\
&\leq& \gamma \left|\frac{\hat d_i - \hat p_i}{a(I-1)}\right| +\gamma \left|\hat\zeta - \tilde\lambda_m\right|  \nonumber \\
&\leq& \gamma\cdot\frac{\sigma}{I-1} + \gamma\cdot\frac{\sigma}{I-1} = \frac{2\gamma\sigma}{I-1} \nonumber
\end{eqnarray}
where the second inequality applies Lipschitz continuity of $\tilde f_i(.)$ and the last inequality exploits the fact that $\hat \zeta \in [\zeta^-, ~\zeta^+]$. 

To finish the proof, we apply the standard definition of convergence. For arbitrary $\epsilon>0$, we can identify integer $I_{\epsilon} \geq \frac{2\gamma\sigma}{\epsilon}+1$, such that for all $I \geq I_\epsilon$, we can make $\left|\hat p_i - \tilde p_i \right| \leq \frac{2\gamma\sigma}{I-1} \leq \epsilon$. This proves $\lim_{I\rightarrow\infty} \left|\hat p_i -\tilde p_i\right|= 0$ for all $i \in \mathcal{I}$. A similar argument can prove  $\lim_{I\rightarrow\infty} \left|\hat d_i -\tilde d_i\right|= 0$ for all $i \in \mathcal{I}$.

\setcounter{equation}{0}  
\renewcommand{\theequation}{D.\arabic{equation}}
\subsection{Proof of Lemma \ref{lemma1}}
\label{apen-4}
The Hessian matrix of $\phi(y)$ is $\mathbb{H}(\phi) = H_1 + H_2 + H_3$,
where
\begin{align}
    H_1 = ~ &  \begin{pmatrix}  \ddot{f}_1 &  & & &  \\ & -\ddot{u}_1 &  & &\\ & & \ddots & & \\ & & & \ddot{f}_I & \\ & & & & -\ddot{u}_I \end{pmatrix} \nonumber\\
    H_2 = ~ & \begin{pmatrix}  \frac{1}{a(I-1)} & \frac{-1}{a(I-1)} & & &  \\ \frac{-1}{a(I-1)} & \frac{1}{a(I-1)} &  & &\\ & & \ddots & & \\ & & & \frac{1}{a(I-1)} & \frac{-1}{a(I-1)}\\ & & & \frac{-1}{a(I-1)} & \frac{1}{a(I-1)} \end{pmatrix} \nonumber\\
    H_3 = ~ & \begin{pmatrix}  \frac{-1}{aI} & \frac{1}{aI} & \ldots & \frac{-1}{aI} & \frac{1}{aI}  \\ \frac{1}{aI} & \frac{-1}{aI} & \ldots & \frac{1}{aI} & \frac{-1}{aI} \\ \vdots & \vdots & & \vdots & \vdots\\   \frac{-1}{aI} & \frac{1}{aI} & \ldots & \frac{-1}{aI} & \frac{1}{aI}  \\ \frac{1}{aI} & \frac{-1}{aI} & \ldots & \frac{1}{aI} & \frac{-1}{aI}  \end{pmatrix} \nonumber
\end{align}

The only non-zero eigenvalue of $H_2$ is $\frac{2}{a(I-1)}$, corresponding to orthonormal eigenvectors $e_i=[0 \cdots 0 ~\underset{(2i-1)}{\frac{\sqrt{2}}{2}}~ \underset{(2i)}{\frac{-\sqrt{2}}{2}}~ 0 \cdots 0]^T,\forall i=1\cdots I$. The only non-zero eigenvalue of $H_3$ is $-\frac{2}{a}$, corresponding to unit eigenvector $e=[\frac{1}{\sqrt{2I}}~\frac{-1}{\sqrt{2I}}\cdots \frac{1}{\sqrt{2I}}~\frac{-1}{\sqrt{2I}}]^T$. When A4 holds, for any vector $x=[x_{11}~x_{12}\cdots x_{i1}~x_{i2} \cdots x_{I1}~x_{I2}]^T \in \mathbb{R}^{2I \times 1}$, we have
\begin{align}
\label{eq:hessian}
    ~ & x^T \mathbb{H}(\phi) x \nonumber\\
    =~ & x^T H_1 x + x^T H_2 x + x^T H_3 x \nonumber\\
    =~ & \sum \limits_{i=1}^I (\ddot{f}_i x_{i1}^2-\ddot{u}_i x_{i2}^2)+  \frac{2}{a(I-1)} (\frac{\sqrt{2}}{2})^2 \sum \limits_{i=1}^T (x_{i1}-x_{i2})^2 \nonumber\\
    ~ & -\frac{2}{a}\frac{1}{2I}(\sum \limits_{i=1}^I \left(x_{i1}-x_{i2})\right)^2 \nonumber\\
    \ge ~ & \sum \limits_{i=1}^I (\ddot{f}_i x_{i1}^2-\ddot{u}_i x_{i2}^2) +  \left(\frac{1}{a(I-1)}-\frac{1}{a}\right) \sum \limits_{i=1}^T (x_{i1}-x_{i2})^2 \nonumber\\
    = ~ & \sum \limits_{i=1}^I \left(\ddot{f}_i x_{i1}^2 - \ddot{u}_i x_{i2}^2 - \frac{I-2}{a(I-1)}(x_{i1}-x_{i2})^2\right) \nonumber\\
    \ge ~& \sum \limits_{i=1}^I \left(\ddot{f}_i x_{i1}^2 - \ddot{u}_i x_{i2}^2 - \frac{2I-4}{a(I-1)}(x_{i1}^2+x_{i2}^2)\right) \ge 0
\end{align}

Therefore, $\mathbb{H}(\phi)$ is a positive semidefinite matrix, implying $\phi(y)$ is a convex function.

Suppose $(\hat{p},\hat{d},\hat{b})$ is an GNE of the game $\mathcal{G}$, and $\hat{\lambda}:=\sum_{i=1}^I \hat{b}_i/(aI)$. According to the KKT condition \eqref{eq:sharing-kkt} and the convexity of $\phi(y)$, it is easy to check that $(\hat{y},\hat{\lambda})$ satisfies
\begin{align}
    L_{\lambda \in \mathbb{R}} (\hat{y},\lambda) \le L(\hat{y},\hat{\lambda}) \le L_{y \in \mathcal{Y}}(y,\hat{\lambda})
\end{align}
which means $(\hat{y},\hat{\lambda})$ is a saddle point of $L(y,\lambda)$.
\end{document}